\documentclass{amsart}
\usepackage{amssymb}

\usepackage{graphicx}

\input{tcilatex}
\input tcilatex

\begin{document}
\title[A birth-and-death approach]{Did the ever dead outnumber the living and when? A birth-and-death approach}
\author{Jean Avan, Nicolas Grosjean and Thierry Huillet$^{*}$}
\address{CNRS, UMR-8089 and University of Cergy-Pontoise\\
2, rue Adolphe Chauvin F-95302, Cergy-Pontoise, Cedex, FRANCE\\
E-mail(s): avan@u-cergy.fr, nicolas.grosjean@u-cergy.fr, huillet@u-cergy.fr}
\maketitle

\begin{abstract}
This paper is an attempt to formalize analytically the question raised in
``World Population Explained: Do Dead People Outnumber Living, Or Vice
Versa?'' Huffington Post, \cite{HJ}. We start developing simple
deterministic Malthusian growth models of the problem (with birth and death
rates either constant or time-dependent) before running into both linear
birth and death Markov chain models and age-structured models.\newline

\textbf{Keywords: }population growth; Malthusian; constant vs time-dependent
birth/death rates; time-inhomogeneous Markov chain; age-structured models.%
\newline

$^{*}$ Corresponding author.
\end{abstract}

\section{Introduction}

Starting from a question raised in (of all places) Huffington Post, asking
whether the total number of dead people ever outgrew the total number of
living people at a certain time, and when the crossover (if any) occurred,
leads us to analyze some simplified models of population growth based on
``mean-field'' approximations, i.e. dealing with the notion of mean
birth/death rates.

We shall successively consider a simple deterministic growth model for the
whole population (Section $2$); a linear birth/death Markov chain model for
the growth of the whole population (Section $3$); a deterministic
age-structured model (Lotka-McKendrick-von Foerster) in Section $4$, for
which we shall explicitly (as far as possible) derive general solutions for
the population sizes (in all 3 cases), the extinction probability, time to
extinction and joint law for the number of people ever born and ever dead
(for the Markovian process). Some explicitly solvable cases will be
considered.

\section{Simple birth and death models for population growth (deterministic)}

We get first interested in the simple deterministic evolution of a
population whose individuals are potentially immortal, non-interacting and
living in a bath with unlimited resources.

\subsection{Constant rates (Malthusian model)}

Let the mean size of some population at time $t=0$ be $x\left( 0\right)
=x_{0}$ and let $x\left( t\right) $ be its mean size at time $t.$

Let $\left( \lambda _{b},\lambda _{d}\right) $ be the birth and death rates
per individual, assumed constant in the first place. With $^{.}=d/dt$,
consider the evolution 
\begin{equation}
\overset{.}{x}\left( t\right) =\left( \lambda _{b}-\lambda _{d}\right)
x\left( t\right) =:\text{ }\overset{.}{x}_{b}\left( t\right) -\overset{.}{x}%
_{d}\left( t\right) ,  \label{rates}
\end{equation}
where $\overset{.}{x}\left( t\right) $ is the growth rate at time $t$ of the
population, whereas $\overset{.}{x}_{b}\left( t\right) =\lambda _{b}x\left(
t\right) $ and $\overset{.}{x}_{d}\left( t\right) =\lambda _{d}x\left(
t\right) $ are the birth and death rates at time $t$ of this population (the
rates at which newborns and dead people are created, respectively).

Integrating, we have $x\left( t\right) =x_{b}\left( t\right) -x_{d}\left(
t\right) $, with $x\left( 0\right) =x_{b}\left( 0\right) .$ The quantities $%
x_{b}\left( t\right) $ and $x_{d}\left( t\right) $ are the mean number of
individuals ever born and ever dead between times $0$ and $t$ and therefore $%
x\left( t\right) $ is the current population size ($x_{b}\left( t\right) $
includes the initial individuals)$.$ Obviously, with $\lambda :=\lambda
_{b}-\lambda _{d}$%
\begin{equation*}
x\left( t\right) =x\left( 0\right) e^{\lambda t}
\end{equation*}
and therefore (if $\lambda _{b}\neq \lambda _{d}$) 
\begin{eqnarray*}
x_{b}\left( t\right) &=&x\left( 0\right) +\frac{\lambda _{b}}{\lambda
_{b}-\lambda _{d}}\left( x\left( t\right) -x\left( 0\right) \right) \\
x_{d}\left( t\right) &=&\frac{\lambda _{d}}{\lambda _{b}-\lambda _{d}}\left(
x\left( t\right) -x\left( 0\right) \right) .
\end{eqnarray*}
Under this model, given that some individual is alive at time $t$, its
lifetime $\tau _{d}$ obeys $\mathbf{P}\left( \tau _{d}>t+\tau \mid \tau
_{d}>t\right) =\mathbf{P}\left( \tau _{d}>\tau \right) =e^{-\lambda _{d}\tau
},$ independent of $t$ (the memory-less property of the exponential
distribution)$.$ The individuals are immortal in that there is no upper
bound for their lifetime although very long lifetimes are unlikely to occur,
i.e. with exponentially small probabilities.

Similarly, its reproduction time $\tau _{b}$ obeys $\mathbf{P}\left( \tau
_{b}>t+\tau \mid \tau _{b}>t\right) =\mathbf{P}\left( \tau _{b}>\tau \right)
=e^{-\lambda _{b}\tau }.$ After each reproduction time, each individual is
replaced by two individuals by binary splitting. Then $\lambda _{d}^{-1}$
and $\lambda _{b}^{-1}$ are the mean lifetime and reproduction times.

We assume the supercritical condition: $\lambda =\lambda _{b}-\lambda _{d}>0$
so that $x\left( t\right) $ grows exponentially fast. Then, as $t\rightarrow
\infty $ 
\begin{equation}
\frac{x_{b}\left( t\right) }{x\left( t\right) }\rightarrow \frac{\lambda _{b}%
}{\lambda _{b}-\lambda _{d}}\text{ and }\frac{x_{d}\left( t\right) }{x\left(
t\right) }\rightarrow \frac{\lambda _{d}}{\lambda _{b}-\lambda _{d}}.
\label{asymp}
\end{equation}

\textbf{Remark:} The cases $\lambda _{b}=\lambda _{d}$ or $\lambda
_{b}-\lambda _{d}<0$ correspond to the critical and subcritical situations,
respectively. In the first case, the population size $x\left( t\right) $
remains constant (while $x_{b}\left( t\right) $ and $x_{d}\left( t\right) $
grow linearly) whereas in the second case, the population size goes to $0$
exponentially fast. Unless otherwise specified we shall stick to the
supercritical case in the sequel.\newline

$\bullet $ Suppose we ask the questions:

- how many people have ever lived in the past and

- do the ever dead outnumber the living and when?

The answer to the first question is $x_{b}\left( t\right) $.

The answer to the second question is positive iff $\lambda _{b}>\lambda _{d}$
and $\frac{\lambda _{d}}{\lambda _{b}-\lambda _{d}}>1$ or else if $\lambda
_{d}<\lambda _{b}<2\lambda _{d}$. Then the time when the overshooting
occurred is $t_{*}$ defined by $\frac{x_{d}\left( t_{*}\right) }{x\left(
t_{*}\right) }=1.$ Thus 
\begin{equation*}
t_{*}=-\frac{1}{\lambda }\log \left( 2-\frac{\lambda _{b}}{\lambda _{d}}%
\right) >0,
\end{equation*}
which is independent of $x\left( 0\right) $. At time $t_{*}$, the number of
people who have ever lived in the past is thus 
\begin{equation*}
x_{b}\left( t_{*}\right) =x\left( 0\right) \left( 1+\frac{\lambda _{b}}{%
\lambda _{b}-\lambda _{d}}\left( e^{\lambda t_{*}}-1\right) \right) =x\left(
0\right) \frac{2\lambda _{d}}{2\lambda _{d}-\lambda _{b}}>x\left( 0\right) .
\end{equation*}
Note also that the population size at $t_{*}$ is 
\begin{equation*}
x\left( t_{*}\right) =x\left( 0\right) e^{\lambda t_{*}}=x\left( 0\right) 
\frac{\lambda _{d}}{2\lambda _{d}-\lambda _{b}}>x\left( 0\right) ,
\end{equation*}
so twice less than $x_{b}\left( t_{*}\right) $ consistently with $x\left(
t_{*}\right) =x_{b}\left( t_{*}\right) -x_{d}\left( t_{*}\right) $ and $%
x_{d}\left( t_{*}\right) =x\left( t_{*}\right) .$

Suppose that at some (large) terminal time $t_{f}$ we know:

- $x_{b}\left( t_{f}\right) $ the number of people who ever lived before
time $t_{f}.$

- the initial population size $x\left( 0\right) $ at time $t=0$.

- the current population size $x\left( t_{f}\right) .$

Then 
\begin{equation*}
\lambda =\lambda _{b}-\lambda _{d}=\frac{1}{t_{f}}\log \left( \frac{x\left(
t_{f}\right) }{x\left( 0\right) }\right) \text{ and }\frac{\lambda _{b}}{%
\lambda _{b}-\lambda _{d}}\sim \frac{x_{b}\left( t_{f}\right) }{x\left(
t_{f}\right) }
\end{equation*}
and both $\lambda _{b}$ and $\lambda _{d}$ are known. As a result, 
\begin{equation*}
t_{*}=-\frac{1}{\lambda }\log \left( 2-\frac{\lambda _{b}}{\lambda _{d}}%
\right) \sim t_{f}\frac{\log \left( 1+x\left( t_{f}\right) /\left(
x_{b}\left( t_{f}\right) -2x\left( t_{f}\right) \right) \right) }{\log
\left( x\left( t_{f}\right) /x\left( 0\right) \right) }
\end{equation*}
can be estimated from the known data.\newline

\textbf{Example:}

The number of people who ever lived on earth is estimated to $x_{b}\left(
t_{f}\right) =105.10^{9}.$ Suppose $10.000$ years ago the initial population
was of $10^{6}$ units. Then $x\left( 0\right) =10^{6}$ and $t_{f}=10^{4}.$
The current world population is around $x\left( t_{f}\right) =7.10^{9}.$
Then 
\begin{equation*}
\lambda _{b}\sim 15.10^{-3}\text{, }\lambda _{d}\sim 14.10^{-3}\text{ and }%
t_{*}\sim 75\text{ years.}
\end{equation*}
These rough estimates based on most simple Malthus growth models are very
sensitive to the initial condition. $\Box $

\subsection{Time-dependent rates}

Because the assumption of constant birth and death rates is not a realistic
hypothesis, we now allow time-dependent rates, so both $\lambda _{b}$ and $%
\lambda _{d}$ are now assumed to depend on current time $t$. We shall assume
that $\lambda _{b}\left( t\right) $ is bounded above by $\lambda
_{b}^{+}<\infty $ and below by $\lambda _{b}^{-}>0$ and similarly for $%
\lambda _{d}\left( t\right) $ with upper and lower bounds $\lambda _{d}^{\pm
}.$ Unless specified otherwise, we assume supercriticality throughout, that
is $\lambda _{b}^{-}>\lambda _{d}^{-}.$ We shall also assume that both $%
\lambda _{b}\left( t\right) $ and $\lambda _{d}\left( t\right) $ are
nonincreasing functions of $t$.\footnote{%
Under this hypothesis, both the expected lifetime and the mean reproduction
time of each individual increase as time passes by, reflecting an
improvement of the general survival conditions.}

Let $\Lambda _{b}\left( t\right) =\int_{0}^{t}ds\lambda _{b}\left( s\right)
, $ $\Lambda _{d}\left( t\right) =\int_{0}^{t}ds\lambda _{d}\left( s\right) $
be primitives of $\lambda _{b}\left( t\right) $ and $\lambda _{d}\left(
t\right) .$

Under this model, given some individual is alive at time $t$, its lifetime $%
\tau _{d}$ obeys $\mathbf{P}\left( \tau _{d}>t+\tau \mid \tau _{d}>t\right)
=e^{-\left( \Lambda _{d}\left( t+\tau \right) -\Lambda _{d}\left( t\right)
\right) }$ so $\mathbf{P}\left( \tau _{d}\in dt\mid \tau _{d}>t\right)
=\lambda _{d}\left( t\right) dt$.

Similarly, its reproduction time $\tau _{b}$ obeys $\mathbf{P}\left( \tau
_{b}>t+\tau \mid \tau _{b}>t\right) =e^{-\left( \Lambda _{b}\left( t+\tau
\right) -\Lambda _{b}\left( t\right) \right) }$ so $\mathbf{P}\left( \tau
_{b}\in dt\mid \tau _{b}>t\right) =\lambda _{b}\left( t\right) dt.$ After
each reproduction time, each individual is again replaced by two
individuals, giving birth to a new individual by binary splitting.

With $^{^{\bullet }}\equiv d/dt$, we thus have to consider the new evolution 
\begin{equation}
\overset{.}{x}\left( t\right) =\left( \lambda _{b}\left( t\right) -\lambda
_{d}\left( t\right) \right) x\left( t\right) =:\text{ }\overset{.}{x}%
_{b}\left( t\right) -\overset{.}{x}_{d}\left( t\right) ,  \label{rates2}
\end{equation}
where $\overset{.}{x}\left( t\right) $ is the growth rate at time $t$ of the
population, whereas $\overset{.}{x}_{b}\left( t\right) =\lambda _{b}\left(
t\right) x\left( t\right) $ and $\overset{.}{x}_{d}\left( t\right) =\lambda
_{d}\left( t\right) x\left( t\right) $ are again the birth and death rates
at time $t$ of this population.

We immediately have $x\left( t\right) =x_{b}\left( t\right) -x_{d}\left(
t\right) $, with $x\left( 0\right) =x_{b}\left( 0\right) .$ The quantities $%
x_{b}\left( t\right) $ and $x_{d}\left( t\right) $ are the mean number of
individuals ever born and ever dead between times $0$ and $t$ ($x_{b}\left(
t\right) $ including the initial individuals) and therefore $x\left(
t\right) $ is the current population size$.$ Clearly, with $\Lambda \left(
t\right) =\Lambda _{b}\left( t\right) -\Lambda _{d}\left( t\right) $%
\begin{equation*}
x\left( t\right) =x\left( 0\right) e^{\Lambda \left( t\right) }
\end{equation*}
and therefore 
\begin{eqnarray*}
x_{b}\left( t\right) &=&x\left( 0\right) \left( 1+\int_{0}^{t}ds\lambda
_{b}\left( s\right) e^{\Lambda \left( s\right) }\right) \\
x_{d}\left( t\right) &=&x\left( 0\right) \int_{0}^{t}ds\lambda _{d}\left(
s\right) e^{\Lambda \left( s\right) }.
\end{eqnarray*}
Note $x\left( t\right) =x\left( 0\right) \left( 1+\int_{0}^{t}ds\left(
\lambda _{b}\left( s\right) -\lambda _{d}\left( s\right) \right) e^{\Lambda
\left( s\right) }\right) $

We assume that $\Lambda \left( t\right) $ is non-decreasing with $t$ so that 
$x\left( t\right) $ grows and also that, with $\lambda _{b}^{-}>\lambda
_{d}^{-}>0$%
\begin{equation*}
\lambda _{b}\left( t\right) \rightarrow \lambda _{b}^{-}>0\text{ and }%
\lambda _{d}\left( t\right) \rightarrow \lambda _{d}^{-}>0\text{ as }%
t\rightarrow \infty .
\end{equation*}
Then, as $t\rightarrow \infty $ 
\begin{equation}
\frac{x_{b}\left( t\right) }{x\left( t\right) }=\frac{1+\int_{0}^{t}ds%
\lambda _{b}\left( s\right) e^{\Lambda \left( s\right) }}{%
1+\int_{0}^{t}ds\left( \lambda _{b}\left( s\right) -\lambda _{d}\left(
s\right) \right) e^{\Lambda \left( s\right) }}\rightarrow \frac{\lambda
_{b}^{-}}{\lambda _{b}^{-}-\lambda _{d}^{-}}\text{ and }  \label{asymp2}
\end{equation}
\begin{equation*}
\frac{x_{d}\left( t\right) }{x\left( t\right) }=\frac{\int_{0}^{t}ds\lambda
_{d}\left( s\right) e^{\Lambda \left( s\right) }}{1+\int_{0}^{t}ds\left(
\lambda _{b}\left( s\right) -\lambda _{d}\left( s\right) \right) e^{\Lambda
\left( s\right) }}\rightarrow \frac{\lambda _{d}^{-}}{\lambda
_{b}^{-}-\lambda _{d}^{-}}.
\end{equation*}
Indeed, by L'Hospital rule, if $\frac{f\left( t\right) }{g\left( t\right) }%
\underset{t\rightarrow \infty }{\rightarrow }c$ $\Rightarrow $ $\frac{%
F\left( t\right) }{G\left( t\right) }\underset{t\rightarrow \infty }{%
\rightarrow }c$ where $\left( F,G\right) $ are the primitives of $\left(
f,g\right) .$ Applying L'Hospital rule to $f=\lambda _{b}\left( t\right)
e^{\Lambda \left( t\right) }$ and $g=\left( \lambda _{b}\left( t\right)
-\lambda _{d}\left( t\right) \right) e^{\Lambda \left( t\right) }$ gives the
result$.$ Suppose $\lambda _{b}^{-}>\lambda _{d}^{-}$. Whenever $\frac{%
\lambda _{d}^{-}}{\lambda _{b}^{-}-\lambda _{d}^{-}}>1$ or else $\lambda
_{b}^{-}<2\lambda _{d}^{-}$, there is a unique $t_{*}$ such that $\frac{%
x_{d}\left( t_{*}\right) }{x\left( t_{*}\right) }=1.$\newline

\textbf{Remark:} Assume $\lambda _{b}^{-}>\lambda _{d}^{-}>0.$ We have the
detailed balance equations: 
\begin{equation*}
\frac{d}{dt}\left[ 
\begin{array}{c}
x_{b}\left( t\right) \\ 
x_{d}\left( t\right)
\end{array}
\right] =\left[ 
\begin{array}{cc}
\lambda _{b}\left( t\right) & -\lambda _{b}\left( t\right) \\ 
\lambda _{d}\left( t\right) & -\lambda _{d}\left( t\right)
\end{array}
\right] \left[ 
\begin{array}{c}
x_{b}\left( t\right) \\ 
x_{d}\left( t\right)
\end{array}
\right] \text{, with }\left[ 
\begin{array}{c}
x_{b}\left( 0\right) \\ 
x_{d}\left( 0\right)
\end{array}
\right] =\left[ 
\begin{array}{c}
x_{b}\left( 0\right) \\ 
0
\end{array}
\right] .
\end{equation*}
In vector form, this is also $\overset{.}{\mathbf{x}}\left( t\right)
=B\left( t\right) \mathbf{x}\left( t\right) $, $\mathbf{x}\left( 0\right) $
where $\mathbf{x}\left( t\right) \succeq \mathbf{0}$ (componentwise)$.$ The
matrix $B\left( t\right) $ has eigenvalues $\left( 0,\lambda _{b}\left(
t\right) -\lambda _{d}\left( t\right) \right) ,$ with $\left( 0,\lambda
_{b}\left( t\right) -\lambda _{d}\left( t\right) \right) \rightarrow \left(
0,\lambda _{b}^{-}-\lambda _{d}^{-}\right) $ as $t\rightarrow \infty .$ $%
\lambda _{1}\left( t\right) =\lambda _{b}\left( t\right) -\lambda _{d}\left(
t\right) $ is its dominant eigenvalue. Letting $X\left( t\right) $ be a $%
2\times 2$ matrix with $X\left( 0\right) =I$, and $\overset{.}{X}\left(
t\right) =B\left( t\right) X\left( t\right) ,$ we have $\mathbf{x}\left(
t\right) =X\left( t\right) \mathbf{x}\left( 0\right) $ and 
\begin{equation*}
e^{-\int_{0}^{t}\lambda _{1}\left( s\right) ds}X\left( t\right) \rightarrow 
\frac{\mathbf{xy}^{^{\prime }}}{\mathbf{x}^{\prime }\mathbf{y}}\text{ as }%
t\rightarrow \infty ,
\end{equation*}
with $B\left( \infty \right) \mathbf{x}=\theta \mathbf{x}$, $\mathbf{y}%
^{\prime }B\left( \infty \right) =\theta \mathbf{y}^{\prime }$, $\theta
=\lambda _{b}^{-}-\lambda _{d}^{-}:=\lambda _{1}\left( \infty \right) $ the
dominant eigenvalue of $B\left( \infty \right) $. The left and right
eigenvectors of $B\left( \infty \right) $ associated to $\theta $ are found
to be $\mathbf{x}^{\prime }=\left( \lambda _{b}^{-},\lambda _{d}^{-}\right) $
and $\mathbf{y}^{\prime }=\left( 1,-1\right) .$ The latter convergence
result holds because, with $\mathbf{y}^{\prime }$ independent of $t$, $%
\mathbf{y}^{\prime }B\left( t\right) =\lambda _{1}\left( t\right) \mathbf{y}%
^{\prime }$ for all $t>0$ and $B\left( t\right) \mathbf{1}=\mathbf{0}$ ($%
\mathbf{1}$ is in the kernel of $B\left( t\right) $ for all $t$), see \cite
{Cohen}. We are thus led to 
\begin{eqnarray*}
e^{-\int_{0}^{t}\lambda _{1}\left( s\right) ds}X\left( t\right)
&=&e^{-\Lambda \left( t\right) }X\left( t\right) \rightarrow \frac{1}{%
\lambda _{b}^{-}-\lambda _{d}^{-}}\left[ 
\begin{array}{cc}
\lambda _{b}^{-} & -\lambda _{b}^{-} \\ 
\lambda _{d}^{-} & -\lambda _{d}^{-}
\end{array}
\right] \text{ as }t\rightarrow \infty \\
\text{and }e^{-\Lambda \left( t\right) }\mathbf{x}\left( t\right)
&\rightarrow &\frac{x\left( 0\right) }{\lambda _{b}^{-}-\lambda _{d}^{-}}%
\left[ 
\begin{array}{c}
\lambda _{b}^{-} \\ 
\lambda _{d}^{-}
\end{array}
\right] \text{ as }t\rightarrow \infty ,
\end{eqnarray*}
which is consistent with the previous analysis making use of L'Hospital
rule, recalling $x\left( t\right) =x_{b}\left( t\right) -x_{d}\left(
t\right) =x\left( 0\right) e^{\Lambda \left( t\right) }$. Note 
\begin{equation*}
\frac{1}{t}\int_{0}^{t}\lambda _{1}\left( s\right) ds\rightarrow \theta 
\text{ as }t\rightarrow \infty .
\end{equation*}

\textbf{Examples:}

$\left( i\right) $ (homographic rates) With $b>a>0$ and $d>c>0,$ assume 
\begin{equation*}
\lambda _{b}\left( t\right) =\frac{at+b}{t+1}\text{ and }\lambda _{d}\left(
t\right) =\frac{ct+d}{t+1}
\end{equation*}
Thus $\lambda _{b}^{-}=a$, $\lambda _{b}^{+}=b$ and $\lambda _{d}^{-}=c$, $%
\lambda _{d}^{+}=d$. Then 
\begin{equation*}
\Lambda _{b}\left( t\right) =at+\left( b-a\right) \log \left( t+1\right) 
\text{ and }\Lambda _{d}\left( t\right) =ct+\left( d-c\right) \log \left(
t+1\right)
\end{equation*}
and with $\delta =\left( b-a\right) -\left( d-c\right) $%
\begin{equation*}
\Lambda \left( t\right) =\left( a-c\right) t+\delta \log \left( t+1\right) 
\text{ and }x\left( t\right) =x\left( 0\right) \left( t+1\right) ^{\delta
}e^{\left( a-c\right) t}.
\end{equation*}
Note $\delta \gtrless 0$, depending on $b\gtrless d$. Thus 
\begin{eqnarray*}
x_{b}\left( t\right) &=&x\left( 0\right) \left( 1+\int_{0}^{t}ds\left(
as+b\right) \left( s+1\right) ^{\delta -1}e^{\left( a-c\right) s}\right) \\
x_{d}\left( t\right) &=&x\left( 0\right) \int_{0}^{t}ds\left( cs+d\right)
\left( s+1\right) ^{\delta -1}e^{\left( a-c\right) s}
\end{eqnarray*}
and $\frac{x_{d}\left( t\right) }{x\left( t\right) }\underset{t\rightarrow
\infty }{\rightarrow }\frac{c}{a-c}.$ When $a>c$, $\log x\left( t\right)
/t\rightarrow a-c>0$ and the population grows at exponential rate.

If in addition $a<2c$, there is a unique $t_{*}$ such that $x_{d}\left(
t_{*}\right) =x\left( t_{*}\right) :$%
\begin{equation*}
\int_{0}^{t_{*}}ds\left( cs+d\right) \left( s+1\right) ^{\delta -1}e^{\left(
a-c\right) s}=\left( t_{*}+1\right) ^{\delta }e^{\left( a-c\right) t_{*}}.
\end{equation*}

$\left( ii\right) $ Gompertz (a critical case $\lambda _{b}^{-}=a=c=\lambda
_{d}^{-}$): We briefly illustrate here on an example that population growth
under criticality also has a rich structure. Let 
\begin{equation*}
\lambda _{b}\left( t\right) =a+\left( b-a\right) e^{-\alpha t}\text{ and }%
\lambda _{d}\left( t\right) =c+\left( d-c\right) e^{-\alpha t}.
\end{equation*}
Assume $\alpha >0,$ $a=c$ and $b-d>\alpha .$ Then 
\begin{equation*}
\Lambda _{b}\left( t\right) =at+\frac{b-a}{\alpha }\left( 1-e^{-\alpha
t}\right) \text{ and }\Lambda _{b}\left( t\right) =ct+\frac{d-c}{\alpha }%
\left( 1-e^{-\alpha t}\right)
\end{equation*}
\begin{equation*}
\Lambda \left( t\right) =\Lambda _{b}\left( t\right) -\Lambda _{d}\left(
t\right) =\frac{b-d}{\alpha }\left( 1-e^{-\alpha t}\right) .
\end{equation*}
Because $\lambda _{b}^{-}=a=c=\lambda _{d}^{-}$, $x_{d}\left( t\right) $ and 
$x\left( t\right) $ will never meet.

We have $x\left( t\right) =x\left( 0\right) e^{\Lambda \left( t\right)
}=x\left( 0\right) e^{\left( b-d\right) \left( 1-e^{-\alpha t}\right)
/\alpha }\rightarrow x\left( 0\right) e^{\left( b-d\right) /\alpha }$ as $%
t\rightarrow \infty :$ due to $b-d>\alpha $, after an initial early time of
exponential growth, the population size stabilizes to a limit. This model
bears some resemblance with the time-homogeneous logistic (or Verhulst)
growth model with alternation of fast and slow growth and an inflection
point between the two regimes. $\Box $

\section{Birth and death Markov chain (stochasticity)}

The deterministic dynamics discussed for $\left( x\left( t\right)
,x_{b}\left( t\right) ,x_{d}\left( t\right) \right) $ arise as the mean
values of some continuous-time stochastic integral-valued birth and death
Markov chain for $\left( N\left( t\right) ,N_{b}\left( t\right) ,N_{d}\left(
t\right) \right) $ which we would like now to discuss; $N_{b}\left( t\right)
,N_{d}\left( t\right) $ are now respectively the number of people who ever
lived and died at time $t$ to the origin, while $N\left( t\right)
=N_{b}\left( t\right) -N_{d}\left( t\right) $ is the number of people
currently alive at $t$. Due to stochasticity, some new effects are expected
to pop in, especially the possibility of extinction, even in the
supercritical regime.

\subsection{The current population size $N\left( t\right) $}

\subsubsection{Probability generating function (pgf) of\textbf{\ }$N\left(
t\right) $}

Let $N\left( t\right) \in \left\{ 0,1,2,...\right\} $ with $N\left( 0\right) 
\overset{d}{\sim }\pi _{0}$ for some given initial probability distribution $%
\pi _{0}$. Define the transition matrix of a (linear) birth and death Markov
chain \cite{KMG} as 
\begin{equation}
\mathbf{P}\left( N\left( t+dt\right) =n+1\mid N\left( t\right) =n\right)
=\lambda _{b}\left( t\right) ndt+o\left( dt\right)  \label{tran1}
\end{equation}
\begin{equation*}
\mathbf{P}\left( N\left( t+dt\right) =n-1\mid N\left( t\right) =n\right)
=\lambda _{d}\left( t\right) ndt+o\left( dt\right) ,
\end{equation*}
with state $\left\{ 0\right\} $ therefore absorbing. Let $\phi _{t}\left(
z\right) =\mathbf{EE}\left( z^{N\left( t\right) }\mid N\left( 0\right)
\right) =:\mathbf{E}\left( z^{N\left( t\right) }\right) $ be the pgf of $%
N\left( t\right) $ averaged over $N\left( 0\right) .$ Let $\phi _{0}\left(
z\right) =\mathbf{E}\left( z^{N\left( 0\right) }\right) $ be the pgf of $%
N\left( 0\right) .$ Then, defining 
\begin{equation*}
g\left( t,z\right) =\left( \lambda _{b}\left( t\right) -\lambda _{d}\left(
t\right) \right) \left( z-1\right) +\lambda _{b}\left( t\right) \left(
z-1\right) ^{2},
\end{equation*}
$\phi _{t}\left( z\right) $ obeys the PDE (see \cite{BR}) 
\begin{equation*}
\partial _{t}\phi _{t}\left( z\right) =g\left( t,z\right) \partial _{z}\phi
_{t}\left( z\right) ,\text{ with }\phi _{t=0}\left( z\right) =\phi
_{0}\left( z\right)
\end{equation*}
or equivalently, solves the non-linear Bernoulli ODE problem 
\begin{equation*}
\overset{.}{\varphi }_{t}\left( z\right) =-g\left( t,\varphi _{t}\left(
z\right) \right) \text{, with }\varphi _{t=0}\left( z\right) =z\text{ and }%
\phi _{t}\left( z\right) =\phi _{0}\left( \varphi _{t}^{-1}\left( z\right)
\right) .
\end{equation*}
Setting $\psi \left( t\right) =e^{\Lambda \left( t\right) }$ and $\eta
\left( t\right) =\int_{0}^{t}ds\lambda _{b}\left( s\right) e^{\Lambda \left(
t\right) -\Lambda \left( s\right) }$, both going to $\infty $ as $%
t\rightarrow \infty $, assuming 
\begin{equation*}
\lambda _{b}\left( t\right) \rightarrow \lambda _{b}^{-}>0\text{ and }%
\lambda _{d}\left( t\right) \rightarrow \lambda _{d}^{-}>0\text{ as }%
t\rightarrow \infty
\end{equation*}
and $\lambda _{b}^{-}\neq \lambda _{d}^{-}$, we easily get the solution 
\begin{equation}
\phi _{t}\left( z\right) =\phi _{0}\left( \frac{1-\left( \eta \left(
t\right) -\psi \left( t\right) \right) \left( z-1\right) }{1-\eta \left(
t\right) \left( z-1\right) }\right) ,  \label{sol1}
\end{equation}
where inside $\phi _{0}\left( \cdot \right) $ we recognize $\varphi
_{t}\left( z\right) $ as an homographic pgf, solution to the ODE problem. We
choose for initial condition a `Bernoulli thinned' geometric random variable
with homographic pgf\footnote{$\mathbf{E}\left( z^{N\left( 0\right) }\right)
:=\phi _{0}\left( z\right) =\phi _{G}\left( \phi _{B}\left( z\right) \right) 
$ is the composition of a geometric pgf $\phi _{G}\left( z\right) =\left(
qz\right) /\left( 1-pz\right) $ ($p+q=1$) with a Bernoulli pgf $\phi
_{B}\left( z\right) =q_{0}+p_{0}z$ ($p_{0}+q_{0}=1$). So $N\left( 0\right) 
\overset{d}{=}\sum_{i=1}^{G}B_{i}$ where the $B_{i}$s are iid Bernoulli,
independent of $G$ with $\mathbf{E}\left( N\left( 0\right) \right) =p_{0}/q$
and Var$\left( N\left( 0\right) \right) =p_{0}\left( q_{0}q+p_{0}p\right)
/q^{2}=q_{0}\mathbf{E}\left( N\left( 0\right) \right) +p\mathbf{E}\left(
N\left( 0\right) \right) ^{2}.$} 
\begin{equation*}
\phi _{0}\left( z\right) =\frac{q\left( q_{0}+p_{0}z\right) }{1-p\left(
q_{0}+p_{0}z\right) }=\frac{q+p_{0}q\left( z-1\right) }{q-p_{0}p\left(
z-1\right) };
\end{equation*}
this considerably simplifies (\ref{sol1}), since the composition of two
homographic pgfs is again an homographic pgf. With this initial condition,
we finally get 
\begin{equation*}
\phi _{t}\left( z\right) =\frac{q-\left( q\eta \left( t\right) -p_{0}q\psi
\left( t\right) \right) \left( z-1\right) }{q-\left( q\eta \left( t\right)
+p_{0}p\psi \left( t\right) \right) \left( z-1\right) }.
\end{equation*}
With $x\left( 0\right) =\phi _{0}^{\prime }\left( 1\right) =p_{0}/q>0$ (with 
$^{\prime }$ the derivative with respect to $z$), averaging over $N\left(
0\right) $, we obtain 
\begin{eqnarray*}
x\left( t\right) &=&\mathbf{E}\left( N\left( t\right) \right) =\mathbf{EE}%
\left( N\left( t\right) \mid N\left( 0\right) \right) =\phi _{t}^{\prime
}\left( 1\right) =\phi _{0}^{\prime }\left( 1\right) \psi \left( t\right)
=x\left( 0\right) e^{\Lambda \left( t\right) } \\
\text{Var}\left( N\left( t\right) \right) &=&\phi _{t}^{^{\prime \prime
}}\left( 1\right) +\phi _{t}^{\prime }\left( 1\right) -\phi _{t}^{\prime
}\left( 1\right) ^{2}=x\left( 0\right) ^{2}\left( 2p-1\right) \psi \left(
t\right) ^{2}+x\left( 0\right) \psi \left( t\right) \left( 1+2\eta \left(
t\right) \right) .
\end{eqnarray*}

\textbf{Remark:} For pgfs of the homographic form (\ref{sol1}), it is easy
to see that, averaging over $N\left( 0\right) $, for $n\geq 1$\footnote{$%
\left[ z^{n}\right] \phi _{t}\left( z\right) $ is the coefficient of $z^{n}$
in the series expansion of $\phi _{t}\left( z\right) .$} 
\begin{eqnarray*}
\mathbf{P}\left( N\left( t\right) =n\right) &=&\mathbf{EP}\left( N\left(
t\right) =n\mid N\left( 0\right) \right) \\
&=&\left[ z^{n}\right] \phi _{t}\left( z\right) =C\left( t\right)
^{n-1}\left( \phi _{t}\left( 0\right) C\left( t\right) +D\left( t\right)
\right) ,
\end{eqnarray*}
where, with $\phi _{t}^{\prime }\left( 0\right) =\mathbf{P}\left( N\left(
t\right) =1\right) =\phi _{0}^{\prime }\left( 0\right) \psi \left( t\right)
\left( 1+\eta \left( t\right) \right) ^{-2}$ 
\begin{equation*}
C\left( t\right) :=1-\frac{\phi _{t}^{\prime }\left( 0\right) }{1-\phi
_{t}\left( 0\right) }\text{ and }D\left( t\right) :=\frac{\phi _{t}^{\prime
}\left( 0\right) }{1-\phi _{t}\left( 0\right) }-\phi _{t}\left( 0\right) .
\end{equation*}

\subsubsection{The extinction probability and time till extinction}

The extinction probability reads 
\begin{equation*}
\mathbf{P}\left( N\left( t\right) =0\right) =\phi _{t}\left( 0\right) =1-%
\frac{p_{0}\psi \left( t\right) }{q+q\eta \left( t\right) +p_{0}p\psi \left(
t\right) }.
\end{equation*}
Because $\mathbf{P}\left( N\left( t\right) =0\right) $ is also $\mathbf{P}%
\left( \tau _{e}\leq t\right) $ where $\tau _{e}$ is the extinction time
(the event $\tau _{e}=\infty $ corresponds to non-extinction), we obtain 
\begin{equation*}
\mathbf{P}\left( \tau _{e}>t\right) =\frac{p_{0}\psi \left( t\right) }{%
q+q\eta \left( t\right) +p_{0}p\psi \left( t\right) }.
\end{equation*}
Whenever $\lambda _{b}^{-}>\lambda _{d}^{-}>0$ (supercritical regime), $\eta
\left( t\right) /\psi \left( t\right) \rightarrow \kappa >1$ (as $%
t\rightarrow \infty $) because 
\begin{equation*}
\eta \left( t\right) /\psi \left( t\right) =\int_{0}^{t}ds\lambda _{b}\left(
s\right) e^{-\Lambda \left( s\right) }>\int_{0}^{t}ds\lambda \left( s\right)
e^{-\Lambda \left( s\right) }=1-e^{-\Lambda \left( t\right) }\rightarrow 1.
\end{equation*}
As it can be checked, were both $\lambda _{b}\left( t\right) $ and $\lambda
_{d}\left( t\right) $ be identified with their limiting values $\lambda
_{b}^{-}$ and $\lambda _{d}^{-}$, then: $\eta \left( t\right) /\psi \left(
t\right) \rightarrow 1/\left( 1-\rho \right) =\lambda _{b}^{-}/\left(
\lambda _{b}^{-}-\lambda _{d}^{-}\right) $ where $\rho :=\lambda
_{d}^{-}/\lambda _{b}^{-}$ ($<1$ in the supercritical regime)$.$

Thus, with $\kappa =1/\left( 1-\rho \right) >1$ 
\begin{equation}
\mathbf{P}\left( \tau _{e}>t\right) \rightarrow \mathbf{P}\left( \tau
_{e}=\infty \right) =\frac{p_{0}}{q\kappa +p_{0}p}\text{ as }t\rightarrow
\infty .  \label{ext}
\end{equation}
The extinction probability reads 
\begin{equation*}
\rho _{e}=\mathbf{P}\left( \tau _{e}<\infty \right) =1-\mathbf{P}\left( \tau
_{e}=\infty \right) =\frac{q\left( \kappa -p_{0}\right) }{q\kappa +p_{0}p}%
=\phi _{0}\left( \rho \right) ,
\end{equation*}
$N\left( t\right) $ vanishes at time $\tau _{e}\mid \tau _{e}<\infty $. Note 
$\rho _{e}\rightarrow 0$ as $x\left( 0\right) \rightarrow \infty $ (or $%
q\rightarrow 0$) and $\phi _{t}\left( z\right) \rightarrow \rho _{e}=\phi
_{0}\left( \rho \right) $ as $t\rightarrow \infty .$

We have 
\begin{equation*}
\tau _{e}=\inf \left( t>0:N_{d}\left( t\right) \geq N_{b}\left( t\right)
\right)
\end{equation*}
the first time the ever dead outnumber the ever living. In the randomized
setup, even in the supercritical regime of mean exponential growth for $%
x\left( t\right) =x_{b}\left( t\right) -x_{d}\left( t\right) $, there is a
``no-luck'' possibility that this event occurs and this corresponds to
global extinction of the population at time $\tau _{e}$.

These results answer the somehow related question on if and when the ever
dead ($N_{d}\left( t\right) $) outnumber the ever living ($N_{b}\left(
t\right) $) in the supercritical linear birth and death Markov chain model;
but so far the question on whether the ever dead ($N_{d}\left( t\right) $)
outnumber the living ($N\left( t\right) =N_{b}\left( t\right) -N_{d}\left(
t\right) $) has not been addressed.

\subsection{The joint law of $N_{b}\left( t\right) $ and $N_{d}\left(
t\right) $}

In order to access to a comparative study between $N_{d}\left( t\right) $
and $N\left( t\right) :=N_{b}\left( t\right) -N_{d}\left( t\right) $, we
need to first compute the joint pgf of $N_{b}\left( t\right) $\textbf{\ }and%
\textbf{\ }$N_{d}\left( t\right) $.

\subsubsection{\textbf{Joint pgf of }$N_{b}\left( t\right) $\textbf{\ and }$%
N_{d}\left( t\right) $}

Define $N_{b}\left( t\right) $ and $N_{d}\left( t\right) $ as the number of
people ever born and ever dead till time $t$. We wish now to compute the
joint law of $N_{b}\left( t\right) $ and $N_{d}\left( t\right) $, keeping in
mind $N\left( t\right) =N_{b}\left( t\right) -N_{d}\left( t\right) $ and for
some $N\left( 0\right) =N_{b}\left( 0\right) .$

With $n=n_{b}-n_{d}$, the joint transition probabilities are obtained as 
\begin{equation}
\mathbf{P}\left( N_{b}\left( t+dt\right) =n_{b}+1,N_{d}\left( t+dt\right)
=n_{d}\mid N\left( t\right) =n\right) =\lambda _{b}\left( t\right)
ndt+o\left( dt\right)  \label{rate2}
\end{equation}

\begin{equation*}
\mathbf{P}\left( N_{b}\left( t+dt\right) =n_{b},N_{d}\left( t+dt\right)
=n_{d}+1\mid N\left( t\right) =n\right) =\lambda _{d}\left( t\right)
ndt+o\left( dt\right)
\end{equation*}
allowing to derive the evolution equation of $\mathbf{EP}\left( N_{b}\left(
t\right) =n_{b},N_{d}\left( t\right) =n_{d}\mid N\left( 0\right) \right) $.

Introducing the joint pgf of $N_{b}\left( t\right) $ and $N_{d}\left(
t\right) $ as $\Phi _{t}\left( z_{b},z_{d}\right) :=\mathbf{E}\left(
z_{b}^{N_{b}\left( t\right) }z_{d}^{N_{d}\left( t\right) }\right) $, we
obtain its time evolution as 
\begin{equation}
\partial _{t}\Phi _{t}=\left( \lambda _{b}\left( t\right) \left(
z_{b}-1\right) +\lambda _{d}\left( t\right) \left( z_{d}-1\right) \right)
\left( z_{b}\partial _{z_{b}}-z_{d}\partial _{z_{d}}\right) \Phi _{t},
\label{jointpgf}
\end{equation}
with initial condition $\Phi _{0}\left( z_{b},z_{d}\right) =\phi _{0}\left(
z_{b}\right) $, noting $N_{d}\left( 0\right) =0.$

\subsubsection{\textbf{Correlations}}

Consistently, we have $x_{b}\left( t\right) =\partial _{z_{b}}\Phi
_{t}\left( 1,1\right) $ and $x_{d}\left( t\right) =\partial _{z_{d}}\Phi
_{t}\left( 1,1\right) .$ Taking the derivative of (\ref{jointpgf}) with
respect to $z_{b}$ and also with respect to $z_{d}$ and evaluating the
results at $z_{b}=z_{d}=1$ we are lead to $\overset{.}{x}_{b}\left( t\right)
=\lambda _{b}\left( t\right) \left( x_{b}\left( t\right) -x_{d}\left(
t\right) \right) =\lambda _{b}\left( t\right) x\left( t\right) $, $%
x_{b}\left( 0\right) =x\left( 0\right) $ and to $\overset{.}{x}_{d}\left(
t\right) =\lambda _{d}\left( t\right) x\left( t\right) $, $x_{d}\left(
0\right) =0.$ With 
\begin{equation*}
\mathbf{E}\left( N_{b}\left( t\right) ^{2}\right) :=\partial
_{z_{b}z_{b}}^{2}\Phi _{t}\left( 1,1\right) ,\text{ }\mathbf{E}\left(
N_{d}\left( t\right) ^{2}\right) :=\partial _{z_{d}z_{d}}^{2}\Phi _{t}\left(
1,1\right) \text{, }\mathbf{E}\left( N_{b}\left( t\right) N_{d}\left(
t\right) \right) :=\partial _{z_{b}z_{d}}^{2}\Phi _{t}\left( 1,1\right)
\end{equation*}
we also have 
\begin{equation*}
\frac{d}{dt}\left[ 
\begin{array}{c}
\mathbf{E}\left( N_{b}\left( t\right) ^{2}\right) \\ 
\mathbf{E}\left( N_{d}\left( t\right) ^{2}\right) \\ 
\mathbf{E}\left( N_{b}\left( t\right) N_{d}\left( t\right) \right)
\end{array}
\right] =\left[ 
\begin{array}{ccc}
2\lambda _{b}\left( t\right) & 0 & -2\lambda _{b}\left( t\right) \\ 
0 & -2\lambda _{d}\left( t\right) & 2\lambda _{d}\left( t\right) \\ 
\lambda _{d}\left( t\right) & -\lambda _{b}\left( t\right) & \lambda
_{b}\left( t\right) -\lambda _{d}\left( t\right)
\end{array}
\right] \left[ 
\begin{array}{c}
\mathbf{E}\left( N_{b}\left( t\right) ^{2}\right) \\ 
\mathbf{E}\left( N_{d}\left( t\right) ^{2}\right) \\ 
\mathbf{E}\left( N_{b}\left( t\right) N_{d}\left( t\right) \right)
\end{array}
\right] \text{, }
\end{equation*}
with initial condition $\left[ 
\begin{array}{c}
\mathbf{E}\left( N_{b}\left( 0\right) ^{2}\right) \\ 
\mathbf{E}\left( N_{d}\left( 0\right) ^{2}\right) \\ 
\mathbf{E}\left( N_{b}\left( 0\right) N_{d}\left( 0\right) \right)
\end{array}
\right] =\left[ 
\begin{array}{c}
\mathbf{E}\left( N\left( 0\right) ^{2}\right) \\ 
0 \\ 
0
\end{array}
\right] .$ In vector form, this is also $\overset{.}{\mathbf{x}}\left(
t\right) =C\left( t\right) \mathbf{x}\left( t\right) $, $\mathbf{x}\left(
0\right) $ where $\mathbf{x}\left( t\right) \succeq \mathbf{0}$
(componentwise)$.$

The matrix $C\left( t\right) $ tends to a limit $C\left( \infty \right) $ as 
$t\rightarrow \infty .$

It has eigenvalues $\left( 0,\lambda _{b}\left( t\right) -\lambda _{d}\left(
t\right) ,2\left( \lambda _{b}\left( t\right) -\lambda _{d}\left( t\right)
\right) \right) ,$ with 
\begin{equation*}
\left( 0,\lambda _{b}\left( t\right) -\lambda _{d}\left( t\right) ,2\left(
\lambda _{b}\left( t\right) -\lambda _{d}\left( t\right) \right) \right)
\rightarrow \left( 0,\lambda _{b}^{-}-\lambda _{d}^{-},2\left( \lambda
_{b}^{-}-\lambda _{d}^{-}\right) \right) ,
\end{equation*}
as $t\rightarrow \infty .$ We identify now $\lambda _{1}\left( t\right)
=2\left( \lambda _{b}\left( t\right) -\lambda _{d}\left( t\right) \right) $
as its dominant eigenvalue (the one for which $\lim_{t\rightarrow \infty }%
\frac{1}{t}\int_{0}^{t}\lambda \left( s\right) ds$ is largest). Letting $%
X\left( t\right) $ be a $3\times 3$ matrix with $X\left( 0\right) =I$ and $%
\overset{.}{X}\left( t\right) =C\left( t\right) X\left( t\right) ,$ we have $%
\mathbf{x}\left( t\right) =X\left( t\right) \mathbf{x}\left( 0\right) $ and 
\begin{equation*}
e^{-\int_{0}^{t}\lambda _{1}\left( s\right) ds}X\left( t\right) \rightarrow 
\frac{\mathbf{xy}^{^{\prime }}}{\mathbf{x}^{\prime }\mathbf{y}}\text{ as }%
t\rightarrow \infty ,
\end{equation*}
with $C\left( \infty \right) \mathbf{x}=\theta \mathbf{x}$, $\mathbf{y}%
^{\prime }C\left( \infty \right) =\theta \mathbf{y}^{\prime }$, $\theta
=2\left( \lambda _{b}^{-}-\lambda _{d}^{-}\right) :=\lambda _{1}\left(
\infty \right) $ the dominant eigenvalue of $C\left( \infty \right) $. The
left and right eigenvectors of $C\left( \infty \right) $ associated to $%
\theta $ are found to be $\mathbf{x}^{\prime }=\left( \left( \lambda
_{b}^{-}\right) ^{2},\left( \lambda _{d}^{-}\right) ^{2},\lambda
_{b}^{-}\lambda _{d}^{-}\right) $ and $\mathbf{y}^{\prime }=\left(
1,1,-2\right) .$ The latter convergence result holds because, with $\mathbf{y%
}^{\prime }$ independent of $t$, $\mathbf{y}^{\prime }C\left( t\right)
=\lambda _{1}\left( t\right) \mathbf{y}^{\prime }$ for all $t>0$ and $%
C\left( t\right) \mathbf{1}=\mathbf{0}$ ($\mathbf{1}$ is in the kernel of $%
C\left( t\right) $ for all $t$).

Indeed, 
\begin{equation*}
X\left( t\right) =\mathbf{x}\left( t\right) \mathbf{y}^{\prime }\text{
solves }\overset{.}{X}\left( t\right) =C\left( t\right) X\left( t\right)
\end{equation*}
and from Theorem $1$ in \cite{levinson}, assuming\footnote{%
The homographic birth and death rates satisfy these conditions. More
generally, the conditions hold for any non increasing rate function $\lambda 
$ of the form $\lambda \left( t\right) =a+\left( b-a\right) \mathbf{P}\left(
T>t\right) $ where $b>a>0$ and $\mathbf{P}\left( T>t\right) $ is the tail
probability distribution of any positive random variable $T.$} 
\begin{equation*}
\int^{\infty }\left| \overset{.}{\lambda }_{b}\left( t\right) \right|
dt<\infty \text{ and }\int^{\infty }\left| \overset{.}{\lambda }_{d}\left(
t\right) \right| dt<\infty ,
\end{equation*}
then $\mathbf{x}\left( t\right) \simeq \mathbf{x}e^{\int_{T}^{t}\lambda
_{1}\left( s\right) ds}$ (for large $t\gg T$)$.$

We are thus led to 
\begin{equation}
e^{-\int_{0}^{t}\lambda _{1}\left( s\right) ds}\text{Cov}\left( N_{b}\left(
t\right) ,N_{d}\left( t\right) \right) \rightarrow C_{b,d}-C_{b}C_{d}>0\text{
as }t\rightarrow \infty  \label{asymcorr}
\end{equation}
where $C_{b,d}=\mathbf{E}\left( N_{b}\left( 0\right) ^{2}\right) \left(
\lambda _{b}^{-}\lambda _{d}^{-}\right) /\left( \lambda _{b}^{-}-\lambda
_{d}^{-}\right) ^{2}$, $C_{b}=\mathbf{E}\left( N_{b}\left( 0\right) \right)
\lambda _{b}^{-}/\left( \lambda _{b}^{-}-\lambda _{d}^{-}\right) $ and $%
C_{d}=\mathbf{E}\left( N_{b}\left( 0\right) \right) \lambda _{d}^{-}/\left(
\lambda _{b}^{-}-\lambda _{d}^{-}\right) $ so $C_{b,d}-C_{b}C_{d}=$Var$%
\left( N_{b}\left( 0\right) \right) \left( \lambda _{b}^{-}\lambda
_{d}^{-}\right) /\left( \lambda _{b}^{-}-\lambda _{d}^{-}\right) ^{2}>0$,
translating an asymptotic positive correlation of $\left( N_{b}\left(
t\right) ,N_{d}\left( t\right) \right) .$ Note 
\begin{equation*}
\frac{1}{t}\int_{0}^{t}\lambda _{1}\left( s\right) ds\rightarrow \theta 
\text{ as }t\rightarrow \infty .
\end{equation*}

\subsubsection{\textbf{Joint pgf revisited}}

Introducing the change of variables from $\left( z_{b},z_{d}\right) $ to $%
\xi \equiv \sqrt{z_{b}/z_{d}}$ and $\zeta \equiv \sqrt{z_{b}z_{d}}$, the
right-hand-side of (\ref{jointpgf}) only contains $\partial _{\xi }.$ The
evolution equation for $\widetilde{\Phi }$ in the new variables indeed reads 
\begin{equation*}
\partial _{t}\widetilde{\Phi }_{t}=\left( \lambda _{b}\left( t\right) \zeta
\xi ^{2}-\left( \lambda _{b}\left( t\right) +\lambda _{d}\left( t\right)
\right) \xi +\lambda _{d}\left( t\right) \zeta \right) \partial _{\xi }%
\widetilde{\Phi }_{t}\text{, }\widetilde{\Phi }_{0}\left( \xi ,\zeta \right)
=\phi _{0}\left( \xi \zeta \right) .
\end{equation*}
$\widetilde{\Phi }_{t}$ is therefore obtained as a function of a suitable
combination of $\xi $ and time (solvable by the method of characteristics),
and separately of $\zeta $ in a way determined by initial conditions.

Fixing $\zeta $ to a certain value, one may then seek for an homographic
ansatz 
\begin{equation}
\widetilde{\Phi }_{t}\left( \xi ,\zeta \right) =\widetilde{\Phi }_{0}\left( 
\frac{a\left( t\right) \xi +b\left( t\right) }{c\left( t\right) \xi +d\left(
t\right) },\zeta \right) =\phi _{0}\left( \frac{a\left( t\right) \xi
+b\left( t\right) }{c\left( t\right) \xi +d\left( t\right) }\zeta \right) ,
\label{solphi}
\end{equation}
with $a\left( 0\right) =d\left( 0\right) =1$ and $b\left( 0\right) =c\left(
0\right) =0.$

Plugging this ansatz into the evolution equation for $\widetilde{\Phi }$ and
denoting respectively 
\begin{equation*}
\alpha \left( t\right) =\lambda _{b}\left( t\right) \zeta \text{, }\beta
\left( t\right) =-\left( \lambda _{b}\left( t\right) +\lambda _{d}\left(
t\right) \right) \text{, }\gamma \left( t\right) =\lambda _{d}\left(
t\right) \zeta ,
\end{equation*}
we need to solve $\overset{.}{\xi }\left( t\right) =\alpha \left( t\right)
\xi \left( t\right) ^{2}+\beta \left( t\right) \xi \left( t\right) +\gamma
\left( t\right) $ with initial condition $\xi \left( 0\right) =\xi $ and the
guess $\xi \left( t\right) =\frac{a\left( t\right) \xi +b\left( t\right) }{%
c\left( t\right) \xi +d\left( t\right) }$. We get 
\begin{eqnarray*}
\overset{.}{a}c-a\overset{.}{c} &=&\left( ad-bc\right) \alpha \\
\overset{.}{b}c-b\overset{.}{c}+\overset{.}{a}d-a\overset{.}{d} &=&\left(
ad-bc\right) \beta \\
\overset{.}{b}d-b\overset{.}{d} &=&\left( ad-bc\right) \gamma .
\end{eqnarray*}
Fixing the $CP^{3}$ gauge choice of any homographic form $\frac{a\xi +b}{%
c\xi +d}$ (invariant under $a,b,c,d\rightarrow \lambda a,\lambda b,\lambda
c,\lambda d$) by fixing $ad-bc=1$ yields 
\begin{equation}
\frac{d}{dt}\left[ 
\begin{array}{cc}
a & c \\ 
b & d
\end{array}
\right] =\left[ 
\begin{array}{cc}
\beta /2 & -\alpha \\ 
\gamma & -\beta /2
\end{array}
\right] \left[ 
\begin{array}{cc}
a & c \\ 
b & d
\end{array}
\right] ,  \label{ordexp}
\end{equation}
to be solved as the ordered exp-integral 
\begin{equation*}
\left[ 
\begin{array}{cc}
a & c \\ 
b & d
\end{array}
\right] \left( t\right) =\overleftarrow{\exp }\int_{0}^{t}\left[ 
\begin{array}{cc}
\beta /2 & -\alpha \\ 
\gamma & -\beta /2
\end{array}
\right] \left( s\right) ds.
\end{equation*}

Although this system can rarely be solved explicitly by quadrature for all $%
t $, this in principle solves $\widetilde{\Phi }_{t}\left( \xi ,\zeta
\right) $ hence $\Phi _{t}\left( z_{b},z_{d}\right) $ when back to the
original variables$.$\newline

\textbf{Remark:} We note from the above differential equation (\ref{ordexp})
and the expression of $\left( \alpha ,\beta ,\gamma \right) $ in the
transition matrix that, as functions of $\zeta $%
\begin{equation*}
\left( a,b/\zeta ,c/\zeta ,d\right) \text{ are functions of }\zeta ^{2},
\end{equation*}
or else that (besides $t$) 
\begin{equation}
\zeta a/c\text{ and }\zeta d/b\text{ are functions of }\zeta ^{2}.
\label{scale}
\end{equation}

\subsubsection{\textbf{Asymptotics}}

In matrix form, (\ref{ordexp}) reads $\overset{.}{X}\left( t\right) =A\left(
t\right) X\left( t\right) $ , $X\left( 0\right) =I$ where $A\left( t\right)
=\left[ 
\begin{array}{cc}
\beta /2 & -\alpha \\ 
\gamma & -\beta /2
\end{array}
\right] $ has zero trace and $X\left( t\right) =\left[ 
\begin{array}{cc}
a & c \\ 
b & d
\end{array}
\right] \left( t\right) $. It follows that det$\left( X\right) $ is a
conserved quantity.

Note that, due to $\left| \zeta \right| \leq 1$, $\det \left( A\left(
t\right) \right) =-\left( \lambda _{b}\left( t\right) +\lambda _{d}\left(
t\right) \right) ^{2}/4+\zeta ^{2}\lambda _{b}\left( t\right) \lambda
_{d}\left( t\right) \leq -\left( \lambda _{b}\left( t\right) -\lambda
_{d}\left( t\right) \right) ^{2}/4<0$: For all $t$, the eigenvalues of $%
A\left( t\right) $ are real with opposite sign $\mu _{\pm }\left( t\right)
=\pm \frac{1}{2}\sqrt{\left( \lambda _{b}\left( t\right) +\lambda _{d}\left(
t\right) \right) ^{2}-4\zeta ^{2}\lambda _{b}\left( t\right) \lambda
_{d}\left( t\right) }$. With 
\begin{equation*}
A_{1}=\left[ 
\begin{array}{cc}
-1/2 & -\zeta \\ 
0 & 1/2
\end{array}
\right] \text{ and }A_{2}=\left[ 
\begin{array}{cc}
-1/2 & 0 \\ 
\zeta & 1/2
\end{array}
\right] ,
\end{equation*}
we have $A\left( t\right) =\lambda _{b}\left( t\right) A_{1}+\lambda
_{d}\left( t\right) A_{2}$ with $\left[ A_{1},A_{2}\right] \neq 0$; In fact $%
A_{1},$ $A_{2}$ generate the full $su(2)$ algebra hence it is consistent to
define 
\begin{equation}
X\left( t\right) =\exp \left( \alpha _{1}\left( t\right) L_{1}+\alpha
_{2}\left( t\right) L_{2}+\alpha _{3}\left( t\right) L_{3}\right) ,
\label{solX}
\end{equation}
for some $\alpha _{i}\left( t\right) $ obeying a Wei-Norman type non-linear
differential equation \cite{hui} and $L_{i}$ Pauli matrices. For recent
further developments, see \cite{OJ}.\newline

With $i=1,2$, let $\mathbf{x}_{i}\left( t\right) =X\left( t\right) \mathbf{x}%
_{i}\left( 0\right) $ with $\mathbf{x}_{1}\left( 0\right) ^{\prime }=\left(
1,0\right) \Rightarrow \mathbf{x}_{1}\left( t\right) ^{\prime }=\left(
a,b\right) $ and $\mathbf{x}_{2}\left( 0\right) ^{\prime }=\left( 0,1\right)
\Rightarrow \mathbf{x}_{2}\left( t\right) ^{\prime }=\left( c,d\right) .$

We set 
\begin{eqnarray*}
A\left( t\right) &=&A\left( \infty \right) +\widetilde{A}\left( t\right) \\
A\left( \infty \right) &=&\lambda _{b}^{-}A_{1}+\lambda _{d}^{-}A_{2} \\
\widetilde{A}\left( t\right) &=&\left( \lambda _{b}\left( t\right) -\lambda
_{b}^{-}\right) A_{1}+\left( \lambda _{d}\left( t\right) -\lambda
_{d}^{-}\right) A_{2}.
\end{eqnarray*}
Of course, the matrix $A\left( \infty \right) $ has two real opposite
eigenvalues 
\begin{equation*}
\mu _{\pm }=\pm \frac{1}{2}\sqrt{\left( \lambda _{b}^{-}+\lambda
_{d}^{-}\right) ^{2}-4\zeta ^{2}\lambda _{b}^{-}\lambda _{d}^{-}}
\end{equation*}
with $\mu _{\pm }\left( t\right) \rightarrow \mu _{\pm }$ as $t\rightarrow
\infty .$

From Theorem $1$ in \cite{levinson}, assuming $\int^{\infty }\left| \overset{%
.}{\lambda }_{b}\left( t\right) \right| dt<\infty $ and $\int^{\infty
}\left| \overset{.}{\lambda }_{d}\left( t\right) \right| dt<\infty ,$ with $%
T $ the (large) time after which the process enters its asymptotic regime,
we have 
\begin{equation}
\mathbf{x}_{i}\left( t\right) \sim \exp \left( \int_{T}^{t}\mu _{+}\left(
s\right) ds\right) \left[ 
\begin{array}{cc}
\phi _{+}\left( 1\right) & \phi _{+}\left( 1\right) \\ 
\phi _{+}\left( 2\right) & \phi _{+}\left( 2\right)
\end{array}
\right] \mathbf{x}_{i}\left( T\right) ,\text{ as }t\rightarrow \infty
\label{levin}
\end{equation}
where $\mathbf{\phi }_{+}$ is the right column eigenvector of $A\left(
\infty \right) $ associated to $\mu _{+}$, i.e. $A\left( \infty \right) 
\mathbf{\phi }_{+}=\mu _{+}\mathbf{\phi }_{+}.$ The row vector $\mathbf{\phi 
}_{+}^{\prime }$ is 
\begin{equation*}
\mathbf{\phi }_{+}^{\prime }:=\left( \phi _{+}\left( 1\right) ,\phi
_{+}\left( 2\right) \right) =\left( 1,\left( \mu _{+}+\left( \lambda
_{b}^{-}+\lambda _{d}^{-}\right) /2\right) /\left( -\zeta \lambda
_{b}^{-}\right) \right) .
\end{equation*}
We conclude that, as $t\rightarrow \infty $, both $\left| \mathbf{x}%
_{i}\left( t\right) \right| \rightarrow \infty .$ Recalling det$\left(
X\right) =a\left( t\right) d\left( t\right) -b\left( t\right) c\left(
t\right) =\left| \mathbf{x}_{1}\left( t\right) \right| \left| \mathbf{x}%
_{2}\left( t\right) \right| \sin \left( \widehat{\mathbf{x}_{1}\left(
t\right) ,\mathbf{x}_{2}\left( t\right) }\right) =1$ for all $t,$ it follows
that the angle $\widehat{\mathbf{x}_{1}\left( t\right) ,\mathbf{x}_{2}\left(
t\right) }$ goes to $0$ as $t\rightarrow \infty $ at rate $2\int_{T}^{t}\mu
_{+}\left( s\right) ds>0:$ the vectors $\mathbf{x}_{1}\left( t\right) ,%
\mathbf{x}_{2}\left( t\right) $ become parallel (or anti-parallel) with the
same time evolution $\sim \exp \left( \int_{T}^{t}\mu _{+}\left( s\right)
ds\right) $ which guarantees that they stay homothetical at large $t\gg T$
(as seen from the asymptotic evaluation (\ref{levin})). Therefore both $%
a\left( t\right) /c\left( t\right) $ and $b\left( t\right) /d\left( t\right) 
$ tend to the same limit (homothetical ratio $\mathbf{x}_{1}/\mathbf{x}_{2}$%
) as $t\rightarrow \infty $, namely $\left( a\left( T\right) +b\left(
T\right) \right) /\left( c\left( T\right) +d\left( T\right) \right) $. This
common limit exists and, from the above scaling argument (\ref{scale}), it
can be written in the form say $h\left( \zeta ^{2}\right) /\zeta $, for some
(unknown) function $h$ encoding the finite size effects till time $T$; it
depends on $\zeta $ and $\zeta ^{2}$. Note that, consistently, both $a\left(
t\right) /b\left( t\right) $ and $c\left( t\right) /d\left( t\right) $ also
have a common limit (angular azimut of both $\mathbf{x}_{1},\mathbf{x}_{2}$)
which is 
\begin{equation*}
\phi _{+}\left( 1\right) /\phi _{+}\left( 2\right) =\left( -\zeta \lambda
_{b}^{-}\right) /\left( \mu _{+}+\left( \lambda _{b}^{-}+\lambda
_{d}^{-}\right) /2\right) =u\left( \zeta ^{2}\right) /\zeta ,
\end{equation*}
where $u\left( \zeta ^{2}\right) =\left( 2\mu _{+}-\left( \lambda
_{b}^{-}+\lambda _{d}^{-}\right) \right) /\left( 2\lambda _{d}^{-}\right) .$

So the limit structure of $X\left( t\right) $ is known and as $t\rightarrow
\infty $%
\begin{equation}
\widetilde{\Phi }_{t}\left( \xi ,\zeta \right) \rightarrow \phi _{0}\left(
h\left( \zeta ^{2}\right) \right) .  \label{solphi2}
\end{equation}
The dependence on $\xi $ has disappeared in the limit. Coming back to the
original variables, recalling $\xi \equiv \sqrt{z_{b}/z_{d}}$ and $\zeta
\equiv \sqrt{z_{b}z_{d}}$, as $t\rightarrow \infty $ yields the limiting
shape of $\Phi _{t}\left( z_{b},z_{d}\right) .$ We get

\begin{equation*}
\Phi _{t}\left( z_{b},z_{d}\right) =\mathbf{E}\left( z_{b}^{N_{b}\left(
t\right) }z_{d}^{N_{d}\left( t\right) }\right) \rightarrow _{t\rightarrow
\infty }\mathbf{E}\left( z_{b}^{N_{b}\left( \infty \right)
}z_{d}^{N_{d}\left( \infty \right) }\right) =\phi _{0}\left( h\left(
z_{b}z_{d}\right) \right) .
\end{equation*}

\subsubsection{\textbf{A particular solvable case}}

It turns out that in several case there is an explicit solution to (\ref
{ordexp}), illustrating the general conclusions just discussed. Let us give
an example.

$\bullet $ \textbf{Explicit solution.} Let $\rho >0$ be a fixed number. When 
$\lambda _{d}\left( t\right) =\rho \lambda _{b}\left( t\right) ,$ (\ref
{ordexp}) can be solved explicitly while observing 
\begin{equation*}
\left[ 
\begin{array}{cc}
\beta /2 & -\alpha \\ 
\gamma & -\beta /2
\end{array}
\right] =\lambda _{b}\left( t\right) \left[ 
\begin{array}{cc}
-\frac{\left( 1+\rho \right) }{2} & -\zeta \\ 
\zeta \rho & \frac{\left( 1+\rho \right) }{2}
\end{array}
\right]
\end{equation*}
and diagonalizing the latter matrix. With $\Lambda _{b}\left( t\right)
=\int_{0}^{t}\lambda _{b}\left( s\right) ds$ and putting 
\begin{equation*}
\theta \left( \zeta ^{2}\right) =\frac{1}{2}\sqrt{\left( 1+\rho \right)
^{2}-4\rho \zeta ^{2}},
\end{equation*}
we indeed get $ad-bc=1$ where 
\begin{eqnarray*}
a &=&-\frac{\left( 1+\rho \right) }{2\theta \left( \zeta ^{2}\right) }\sinh
\left( \theta \left( \zeta ^{2}\right) \Lambda _{b}\left( t\right) \right)
+\cosh \left( \theta \left( \zeta ^{2}\right) \Lambda _{b}\left( t\right)
\right) \\
b &=&\frac{\zeta \rho }{\theta \left( \zeta ^{2}\right) }\sinh \left( \theta
\left( \zeta ^{2}\right) \Lambda _{b}\left( t\right) \right) \\
c &=&-\frac{\zeta }{\theta \left( \zeta ^{2}\right) }\sinh \left( \theta
\left( \zeta ^{2}\right) \Lambda _{b}\left( t\right) \right) \\
d &=&\frac{\left( 1+\rho \right) }{2\theta \left( \zeta ^{2}\right) }\sinh
\left( \theta \left( \zeta ^{2}\right) \Lambda _{b}\left( t\right) \right)
+\cosh \left( \theta \left( \zeta ^{2}\right) \Lambda _{b}\left( t\right)
\right) .
\end{eqnarray*}
so that 
\begin{equation*}
\widetilde{\Phi }_{t}\left( \xi ,\zeta \right) =\widetilde{\Phi }_{0}\left( 
\frac{a\xi +b}{c\xi +d},\zeta \right) =\phi _{0}\left( \frac{a\xi +b}{c\xi +d%
}\zeta \right) .
\end{equation*}
Recalling $\xi \zeta =z_{b}$ and $\zeta ^{2}=z_{b}z_{d}$, coming back to the
original variables, we get 
\begin{equation*}
\Phi _{t}\left( z_{b},z_{d}\right) =\mathbf{E}\left( z_{b}^{N_{b}\left(
t\right) }z_{d}^{N_{d}\left( t\right) }\right) =\phi _{0}\left( \frac{%
A_{t}\left( z_{b}z_{d}\right) z_{b}+B_{t}\left( z_{b}z_{d}\right) z_{b}z_{d}%
}{C_{t}\left( z_{b}z_{d}\right) z_{b}+D_{t}\left( z_{b}z_{d}\right) }\right)
\end{equation*}
where 
\begin{eqnarray*}
A_{t}\left( z_{b}z_{d}\right) &=&-\frac{1+\rho }{2\theta \left(
z_{b}z_{d}\right) }\sinh \left( \theta \left( z_{b}z_{d}\right) \Lambda
_{b}\left( t\right) \right) +\cosh \left( \theta \left( z_{b}z_{d}\right)
\Lambda _{b}\left( t\right) \right) \\
B_{t}\left( z_{b}z_{d}\right) &=&\frac{\rho }{\theta \left(
z_{b}z_{d}\right) }\sinh \left( \theta \left( z_{b}z_{d}\right) \Lambda
_{b}\left( t\right) \right) \\
C_{t}\left( z_{b}z_{d}\right) &=&-\frac{1}{\theta \left( z_{b}z_{d}\right) }%
\sinh \left( \theta \left( z_{b}z_{d}\right) \Lambda _{b}\left( t\right)
\right) \\
D_{t}\left( z_{b}z_{d}\right) &=&\frac{1+\rho }{2\theta \left(
z_{b}z_{d}\right) }\sinh \left( \theta \left( z_{b}z_{d}\right) \Lambda
_{b}\left( t\right) \right) +\cosh \left( \theta \left( z_{b}z_{d}\right)
\Lambda _{b}\left( t\right) \right)
\end{eqnarray*}
each depend on $t$ and $z_{b}z_{d}.$

$\bullet $ $t\rightarrow \infty $ \textbf{asymptotics}.

This fully solvable case enables us to obtain exact asymptotics for the
generating function $\Phi _{t}\left( z_{b},z_{d}\right) $ in the
supercritical case $\rho <1$. From the explicit formulae for $\left(
a,b,c,d\right) $ we get:

\begin{equation*}
\lim_{t\rightarrow \infty }\frac{a\left( t\right) }{c\left( t\right) }=\frac{%
1+\rho -2\theta \left( \zeta ^{2}\right) }{2\zeta }=\frac{2\zeta \rho }{%
1+\rho +2\theta \left( \zeta ^{2}\right) }=\lim_{t\rightarrow \infty }\frac{%
b\left( t\right) }{d\left( t\right) }
\end{equation*}
identified indeed by definition of $\theta \left( \zeta ^{2}\right) .$ This
illustrates the fact that the vectors $\mathbf{x}_{i}\left( t\right) $, $%
i=1,2$ become parallel up to $\exp \left( -2\Lambda _{b}\left( t\right)
\right) $ when $t$ gets large. Note that one still have of course $ad-bc=1$ $%
\Leftrightarrow $ $\left( ad\right) /\left( bc\right) -1=1/\left( bc\right)
=O\left( \exp \left( -2\Lambda _{b}\left( t\right) \right) \right) $.
Indeed, it appears that $a/c-b/d=1/\left( cd\right) =O\left( \exp \left(
-2\Lambda _{b}\left( t\right) \right) \right) \rightarrow 0$ (without being
strictly $=0$). Hence we get that the generating function yields 
\begin{equation*}
\Phi _{t}\left( z_{b},z_{d}\right) \rightarrow \phi _{0}\left( \frac{1+\rho
-2\theta }{2}\right) ,\text{ }\theta =\frac{1}{2}\sqrt{\left( 1+\rho \right)
^{2}-4\rho z_{b}z_{d}}.
\end{equation*}
Dependence on $\xi $ has disappeared. Observing that as $t\rightarrow \infty 
$%
\begin{equation}
\Phi _{t}\left( z_{b},z_{d}\right) =\mathbf{E}\left( z_{b}^{N_{b}\left(
t\right) }z_{d}^{N_{d}\left( t\right) }\right) \rightarrow \rho _{e}\mathbf{E%
}\left( \left( z_{b}^{{}}z_{d}\right) ^{N_{d}\left( \infty \right) }\mid 
\text{extinct}\right) ,  \label{asyphi}
\end{equation}
because either extinction occurred (with probability $\rho _{e}$) and $%
N_{b}\left( \infty \right) =N_{d}\left( \infty \right) <\infty $ or not in
which case $N_{b}\left( \infty \right) =N_{d}\left( \infty \right) =\infty ,$
we conclude 
\begin{equation}
\rho _{e}\mathbf{E}\left( z^{N_{d}\left( \infty \right) }\mid \text{extinct}%
\right) =\phi _{0}\left( \frac{1+\rho -\sqrt{\left( 1+\rho \right)
^{2}-4\rho z}}{2}\right) .  \label{asyphi2}
\end{equation}
Recalling $\rho _{e}=\phi _{0}\left( \rho \right) $, we obtain the pgf of
the cumulative population size (ever born or ever dead) $N_{b}\left( \infty
\right) =N_{d}\left( \infty \right) <\infty $, given extinction occurred, as
the compound expression

\begin{eqnarray*}
\mathbf{E}\left( z^{N_{d}\left( \infty \right) }\mid \text{extinct}\right)
&=&\phi _{0}\left( h\left( z\right) \right) /\phi _{0}\left( h\left(
1\right) \right) \\
\text{where }h\left( z\right) &=&\frac{1+\rho -\sqrt{\left( 1+\rho \right)
^{2}-4\rho z}}{2}
\end{eqnarray*}
In particular (with $h\left( 1\right) =\rho =\lambda _{d}^{-}/\lambda
_{b}^{-}<1$) 
\begin{equation*}
\mathbf{E}\left( N_{d}\left( \infty \right) \mid \text{extinct}\right) =%
\frac{\phi _{0}^{\prime }\left( h\left( 1\right) \right) }{\phi _{0}\left(
h\left( 1\right) \right) }h^{\prime }\left( 1\right) =\frac{\phi
_{0}^{\prime }\left( \rho \right) }{\phi _{0}\left( \rho \right) }\frac{\rho 
}{1-\rho }.
\end{equation*}
Note that, would one consider a subcritical case $\rho >1$, then $h\left(
1\right) =1$ leading to $\rho _{e}=\phi _{0}\left( h\left( 1\right) \right)
=1$ (almost sure extinction) with $\mathbf{E}\left( z^{N_{d}\left( \infty
\right) }\right) =\phi _{0}\left( h\left( z\right) \right) $ characterizing
the law of $N_{b}\left( \infty \right) =N_{d}\left( \infty \right) <\infty $
at extinction. Note $\mathbf{E}\left( N_{d}\left( \infty \right) \right) =%
\mathbf{E}\left( N_{b}\left( 0\right) \right) \frac{\rho }{\rho -1}.$

In the critical case $\rho =1$, $h\left( 1\right) =1$ $\left( \rho
_{e}=1\right) $ and $\mathbf{E}\left( z^{N_{d}\left( \infty \right) }\right)
=\phi _{0}\left( 1-\sqrt{1-z}\right) $ with $\mathbf{E}\left( N_{d}\left(
\infty \right) \right) =\infty .$\newline

\textbf{Remark:} This solvable case includes the situation where both $%
\lambda _{b}\left( t\right) $ and $\lambda _{d}\left( t\right) $ are
independent of time. It suffices to particularize the latter formula while
setting $\Lambda _{b}\left( t\right) =\lambda _{b}t.$ $\Box $

If $z_{b}=z,$ $z_{d}=z^{-1},$ 
\begin{equation*}
\Phi _{t}\left( z,z^{-1}\right) =\mathbf{E}\left( z^{N\left( t\right)
}\right) =\phi _{0}\left( \frac{A_{t}\left( 1\right) z+B_{t}\left( 1\right) 
}{C_{t}\left( 1\right) z+D_{t}\left( 1\right) }\right)
\end{equation*}
which is the previous homographic expression for $\mathbf{E}\left(
z^{N\left( t\right) }\right) $ in the particular case when $\lambda
_{d}\left( t\right) =\rho \lambda _{b}\left( t\right) $.

If $z_{b}=z,$ $z_{d}=1$ or $z_{b}=1,$ $z_{d}=z$, we get the marginals 
\begin{eqnarray*}
\Phi _{t}\left( z,1\right) &=&\mathbf{E}\left( z^{N_{b}\left( t\right)
}\right) =\phi _{0}\left( z\frac{A_{t}\left( z\right) +B_{t}\left( z\right) 
}{C_{t}\left( z\right) z+D_{t}\left( z\right) }\right) \text{ or} \\
\Phi _{t}\left( 1,z\right) &=&\mathbf{E}\left( z^{N_{d}\left( t\right)
}\right) =\phi _{0}\left( \frac{A_{t}\left( z\right) +B_{t}\left( z\right) z%
}{C_{t}\left( z\right) +D_{t}\left( z\right) }\right) .
\end{eqnarray*}
The processes $\left( N_{b}\left( t\right) ,N_{d}\left( t\right) \right) $
both are monotone non-decreasing.

If $z_{b}=z^{-1},$ $z_{d}=z^{2}$, we are interested in the Laurent series 
\begin{equation}
\psi _{t}\left( z\right) :=\mathbf{E}\left( z^{2N_{d}\left( t\right)
-N_{b}\left( t\right) }\right) =\Phi _{t}\left( z^{-1},z^{2}\right) =\phi
_{0}\left( \frac{A_{t}\left( z\right) z^{-1}+B_{t}\left( z\right) z}{%
C_{t}\left( z\right) z^{-1}+D_{t}\left( z\right) }\right) .  \label{phit}
\end{equation}

The quantity $\psi _{t}\left( z\right) =\Phi _{t}\left( z^{-1},z^{2}\right) $
turns out to be useful in the understanding of the statistical properties of
the first time that $\Delta \left( t\right) :=2N_{d}\left( t\right)
-N_{b}\left( t\right) $ hits $0$, starting from $\Delta \left( 0\right)
=-N_{b}\left( 0\right) =-N\left( 0\right) <0$. We now address this point.

\subsubsection{\textbf{First time }$N_{d}\left( t\right) \geq N\left(
t\right) =N_{b}\left( t\right) -N_{d}\left( t\right) $\textbf{\ (}$\Delta
\left( t\right) \geq 0$\textbf{)}}

The process $\Delta \left( t\right) $ takes values in $\Bbb{Z}$. It moves up
by 2 units when a death event occurs and down by 1 unit when a birth event
occurs. Starting from $\Delta \left( 0\right) =-N_{b}\left( 0\right) <0,$
the process $\Delta \left( t\right) $ will enter the nonnegative region $%
\Delta \left( t\right) \geq 0$ for the first time at some random crossing
time $\tau _{\Delta \left( 0\right) }$ where the dependence on the initial
condition has been emphasized.

There are two types of crossing events:

- type-1: $\Delta \left( t\right) $ enters the nonnegative region $\Delta
\left( t\right) \geq 0$ as a result of $\Delta \left( t\right) =-1$ followed
instantaneously by a death event leading to $\Delta \left( t_{+}\right) =+1.$

- type-2: $\Delta \left( t\right) $ enters the nonnegative region $\Delta
\left( t\right) \geq 0$ as a result of $\Delta \left( t\right) =-2$ followed
instantaneously by a death event leading to $\Delta \left( t_{+}\right) =0.$

Consider the process $\Delta _{1}\left( t\right) $ which is $\Delta \left(
t\right) $ conditioned on $t=0$ being a type-1 first crossing time ($\Delta
\left( 0\right) =-1$ and $\Delta \left( 0_{+}\right) =+1$). Let $\lambda
_{1}\left( t\right) $ be the rate at which some type-1 crossing occurs at $t$
for $\Delta _{1}.$ We let 
\begin{equation*}
\widehat{\lambda }_{1}\left( s\right) :=\int_{0}^{\infty }e^{-st}\lambda
_{1}\left( t\right) dt
\end{equation*}
be the Laplace-Stieltjes transform (LST) of $\lambda _{1}\left( t\right) .$

Let $\tau _{1}$ denote the (random) time between 2 consecutive type-1
crossing times for $\Delta _{1}.$ By standard renewal arguments \cite{Feller}
\begin{equation}
\widehat{\lambda }_{1}\left( s\right) =1/\left( 1-\widehat{f}_{1}\left(
s\right) \right) ,  \label{renew}
\end{equation}
where $\widehat{f}_{1}\left( s\right) $ is the LST of the law of $\tau _{1}.$

Similarly, consider the process $\Delta _{2}\left( t\right) $ which is $%
\Delta \left( t\right) $ conditioned on $t=0$ being a type-2 first crossing
time ($\Delta \left( 0\right) =-2$ and $\Delta \left( 0_{+}\right) =0$). Let 
$\lambda _{2}\left( t\right) $ be the rate at which some type-2 crossing
occurs at $t$ for $\Delta _{2}.$ Let $\tau _{2}$ denote the random time
between 2 consecutive type-2 crossing times for $\Delta _{2}.$ Then $%
\widehat{\lambda }_{2}\left( s\right) =1/\left( 1-\widehat{f}_{2}\left(
s\right) \right) $ where $\widehat{\lambda }_{1}\left( s\right) $ is the LST
of $\lambda _{2}\left( t\right) $ and $\widehat{f}_{2}\left( s\right) $ the
LST of the law of $\tau _{2}.$

Clearly, between any two consecutive crossings of any type, there can be
none, one or more crossings of the other type.

The laws of the first return times $\left( \tau _{1},\tau _{2}\right) $ of $%
\left( \Delta _{1},\Delta _{2}\right) $ are thus known once the rates $%
\left( \lambda _{1}\left( t\right) ,\lambda _{2}\left( t\right) \right) $
are.\newline

Let now $\lambda _{1,\Delta \left( 0\right) }\left( t\right) $ (respectively 
$\lambda _{2,\Delta \left( 0\right) }\left( t\right) $) be the rate at which
some type-1 (respectively type-2) crossing occurs at $t$ when $\Delta $ is
now simply conditioned on any $\Delta \left( 0\right) <0.$

Let also $\tau _{1,\Delta \left( 0\right) }$ (respectively $\tau _{2,\Delta
\left( 0\right) }$) be the first time that some type-1 (respectively type-2)
crossing occurs for $\Delta $ given $\Delta \left( 0\right) <0.$

Again by renewal arguments for this now delayed renewal process: 
\begin{equation}
\widehat{\lambda }_{1,\Delta \left( 0\right) }\left( s\right) =\widehat{f}%
_{1,\Delta \left( 0\right) }\left( s\right) /\left( 1-\widehat{f}_{1}\left(
s\right) \right) \text{ and}  \label{lst}
\end{equation}
\begin{equation*}
\widehat{\lambda }_{2,\Delta \left( 0\right) }\left( s\right) =\widehat{f}%
_{2,\Delta \left( 0\right) }\left( s\right) /\left( 1-\widehat{f}_{2}\left(
s\right) \right)
\end{equation*}
where $\widehat{\lambda }_{1,\Delta \left( 0\right) }\left( s\right) $
(respectively $\widehat{\lambda }_{2,\Delta \left( 0\right) }\left( s\right) 
$) is the LST of $\lambda _{1,\Delta \left( 0\right) }\left( t\right) $
(respectively $\lambda _{2,\Delta \left( 0\right) }\left( t\right) $) and $%
\widehat{f}_{1,\Delta \left( 0\right) }\left( s\right) $ (respectively $%
\widehat{f}_{2,\Delta \left( 0\right) }\left( s\right) $) is the LST of the
law of $\tau _{1,\Delta \left( 0\right) }$ (respectively $\tau _{2,\Delta
\left( 0\right) }$). Finally, we clearly have 
\begin{equation}
\tau _{\Delta \left( 0\right) }=\inf \left( \tau _{1,\Delta \left( 0\right)
},\tau _{2,\Delta \left( 0\right) }\right) .  \label{firsthit}
\end{equation}

The knowledge of all these quantities (in particular the law of $\tau
_{\Delta \left( 0\right) }$) is therefore conditioned on the knowledge of
the rates. Conditioned on the initial condition $\Delta \left( 0\right) <0$,
we have for instance 
\begin{equation}
\lambda _{1,\Delta \left( 0\right) }\left( t\right) dt=\left[ z^{-1}\right]
\psi _{t}\left( z\right) \cdot \lambda _{d}\left( t\right) \cdot \mathbf{E}%
_{\Delta \left( 0\right) }\left( N_{d}\left( t\right) +1\mid \Delta \left(
t\right) =-1\right) dt  \label{ratesol}
\end{equation}
\begin{equation*}
\lambda _{2,\Delta \left( 0\right) }\left( t\right) dt=\left[ z^{-2}\right]
\psi _{t}\left( z\right) \cdot \lambda _{d}\left( t\right) \cdot \mathbf{E}%
_{\Delta \left( 0\right) }\left( N_{d}\left( t\right) +2\mid \Delta \left(
t\right) =-2\right) dt,
\end{equation*}
corresponding respectively, through a jump of size $2$ of $\Delta \left(
t\right) $, to the death events $\Delta \left( t+dt\right) =1$ given $\Delta
\left( t\right) =-1$ and $\Delta \left( t+dt\right) =0$ given $\Delta \left(
t\right) =-2.$ In (\ref{ratesol}), $\psi _{t}\left( z\right) =\Phi
_{t}\left( z^{-1},z^{2}\right) $ is given by (\ref{phit}).

The term $\left[ z^{-1}\right] \psi _{t}\left( z\right) $ (respectively $%
\left[ z^{-2}\right] \psi _{t}\left( z\right) $) is the probability that $%
\Delta \left( t\right) =-1$ (respectively $\Delta \left( t\right) =-2$).

The term $\mathbf{E}\left( N_{d}\left( t\right) +1\mid \Delta \left(
t\right) =-1\right) $ (respectively $\mathbf{E}\left( N_{d}\left( t\right)
+2\mid \Delta \left( t\right) =-2\right) $) is the expected value of $%
N\left( t\right) =N_{b}\left( t\right) -N_{d}\left( t\right) $ given $\Delta
\left( t\right) =-1$ (respectively of $N\left( t\right) $ given $\Delta
\left( t\right) =-2$).

These last two informations can be extracted from the joint probability
generating function of $\left( N_{d}\left( t\right) ,\Delta \left( t\right)
:=2N_{d}\left( t\right) -N_{b}\left( t\right) \mid \Delta \left( 0\right)
\right) $ which is known from $\Phi _{t}\left( z_{b},z_{d}\right) =\mathbf{E}%
\left( z_{b}^{N_{b}\left( t\right) }z_{d}^{N_{d}\left( t\right) }\right) .$

Therefore, the law of $\tau _{\Delta \left( 0\right) }$ follows in
principle, although the actual computational task left remains huge.

\section{Age-structured models}

We finally address a different point of view of the birth/death problems,
once again deterministic but now based on age-structured models (see e.g. 
\cite{CharlesW}, \cite{Cushing}, \cite{HoppenS} and \cite{Webb}).\newline

Let $\rho =\rho \left( a,t\right) $ be the density of individuals of age $a$
at time $t$. Then $\rho $ obeys 
\begin{equation}
\partial _{t}\rho +\partial _{a}\rho =-\lambda _{d}\left( a,t\right) \rho ,
\label{eq1}
\end{equation}
with boundary conditions 
\begin{equation*}
\rho \left( a,0\right) =\rho _{0}\left( a\right) \text{ and }\rho \left(
0,t\right) =\int \lambda _{b}\left( a,t\right) \rho \left( a,t\right) da.
\end{equation*}
The quantities $\lambda _{b}\left( a,t\right) $ and $\lambda _{d}\left(
a,t\right) $ are the birth and death rates at age and time $\left(
a,t\right) $ and they are the inputs of the model together with $\rho
_{0}\left( a\right) .$ The equation (\ref{eq1}) with its boundary conditions
constitute the Lotka-McKendrick-Von Foerster model, \cite{McK}.

The total population mean size at $t$ is 
\begin{equation*}
x\left( t\right) =\int \rho \left( a,t\right) da.
\end{equation*}
With $^{^{\bullet }}\equiv d/dt$, integrating (\ref{eq1}) with respect to
age and observing $\rho \left( \infty ,t\right) =0$, we get 
\begin{equation*}
\overset{.}{x}\left( t\right) =\int \lambda _{b}\left( a,t\right) \rho
\left( a,t\right) da-\int \lambda _{d}\left( a,t\right) \rho \left(
a,t\right) da=:\text{ }\overset{.}{x}_{b}\left( t\right) -\overset{.}{x}%
_{d}\left( t\right)
\end{equation*}
where $\overset{.}{x}_{b}\left( t\right) $ and $\overset{.}{x}_{d}\left(
t\right) $ are the birth and death rates at time $t$. We note that $\overset{%
.}{x}_{b}\left( t\right) =\rho \left( 0,t\right) $, the rate at which new
individuals (of age $a=0$) are injected in the system$.$

Integrating, we have $x\left( t\right) =x_{b}\left( t\right) -x_{d}\left(
t\right) $, with $x\left( 0\right) =x_{b}\left( 0\right) $ in view of $%
\lambda _{d}\left( a,0\right) =0.$ $x_{b}\left( t\right) $ and $x_{d}\left(
t\right) $ are the mean number of individuals ever born and ever dead
between times $0$ and $t$ ($x_{b}\left( t\right) $ including, as before, the
initial individuals)$.$ The answer to the question of whether and when the
ever dead outnumbered the living is transferred to the existence of the time 
$t_{*}$ when some overshooting occurred, defined by $\frac{x_{d}\left(
t_{*}\right) }{x\left( t_{*}\right) }=1.$ As in Section $2$, this requires
the computation of both $x_{d}\left( t\right) $ and $x\left( t\right) $ and
we now address this question.\newline

With $a\wedge t=\min \left\{ a,t\right\} $, let 
\begin{equation}
\overline{\Lambda }_{d}\left( a,t\right) =\int_{0}^{a\wedge t}\lambda
_{d}\left( a-s,t-s\right) ds.  \label{eq2.0}
\end{equation}
Integration of (\ref{eq1}), using the method of characteristics yields its
solution as 
\begin{equation}
\rho \left( a,t\right) =\rho _{0}\left( a-t\right) e^{-\overline{\Lambda }%
^{d}\left( a,t\right) }\text{ if }t<a  \label{eq2}
\end{equation}
\begin{equation*}
=\text{ }\overset{.}{x}_{b}\left( t-a\right) e^{-\overline{\Lambda }%
^{d}\left( a,t\right) }\text{ if }t>a.
\end{equation*}
Recall $\rho _{0}\left( a\right) =\rho \left( a,0\right) $ and $\overset{.}{x%
}_{b}\left( t\right) =\rho \left( 0,t\right) $ as from the boundary
conditions of (\ref{eq1})$.$ We shall distinguish 3 cases.

\subsection{Case: $\lambda _{b}\left( a,t\right) $ and $\lambda _{d}\left(
a,t\right) $ independent of age $a$\protect\footnote{%
Not very realistic since it assumes that individuals may reproduce at any
age and are potentially immortal eg no cutoff to $\lambda _{b}\left(
a\right) $ for $a<a_{crit}$ or $a>a_{crit}$, and $\lambda _{d}\left(
a\right) $ does not increase with $a.$}}

Suppose $\lambda _{b}\left( a,t\right) =\lambda _{b}\left( t\right) $ and $%
\lambda _{d}\left( a,t\right) =\lambda _{d}\left( t\right) $ for some given
time-dependent rate functions $\lambda _{b}\left( t\right) ,\lambda
_{d}\left( t\right) .$ Then $\overset{.}{x}_{b}\left( t\right) =\lambda
_{b}\left( t\right) x\left( t\right) $ and $\overset{.}{x}_{d}\left(
t\right) =\lambda _{d}\left( t\right) x\left( t\right) .$ Supposing $x\left(
0\right) =\int \rho _{0}\left( a\right) da<\infty $, we get 
\begin{equation*}
x\left( t\right) =x\left( 0\right) e^{\Lambda _{b}\left( t\right) -\Lambda
_{d}\left( t\right) },
\end{equation*}
where $\Lambda _{b},\Lambda _{d}$ are the primitives of $\lambda
_{b},\lambda _{d}.$ Furthermore, 
\begin{eqnarray*}
\overline{\Lambda }_{d}\left( a,t\right) &=&\int_{0}^{a\wedge t}\lambda
_{d}\left( t-s\right) ds \\
&=&\Lambda _{d}\left( t\right) \text{ if }t<a \\
&=&\Lambda _{d}\left( t\right) -\Lambda _{d}\left( t-a\right) \text{ if }t>a.
\end{eqnarray*}
We thus obtain 
\begin{eqnarray*}
x_{b}\left( t\right) &=&x\left( 0\right) +\int_{0}^{t}\lambda _{b}\left(
s\right) x\left( s\right) ds=x\left( 0\right) \left( 1+\int_{0}^{t}\lambda
_{b}\left( s\right) e^{\Lambda _{b}\left( s\right) -\Lambda _{d}\left(
s\right) }ds\right) \\
x_{d}\left( t\right) &=&\int_{0}^{t}\lambda _{d}\left( s\right) x\left(
s\right) ds=x\left( 0\right) \int_{0}^{t}\lambda _{d}\left( s\right)
e^{\Lambda _{b}\left( s\right) -\Lambda _{d}\left( s\right) }ds
\end{eqnarray*}
and for the age-structure 
\begin{eqnarray*}
\rho \left( a,t\right) &=&\rho _{0}\left( a-t\right) e^{-\Lambda _{d}\left(
t\right) }\text{ if }t<a \\
&=&\overset{.}{x}_{b}\left( t-a\right) e^{-\left( \Lambda _{d}\left(
t\right) -\Lambda _{b}\left( t-a\right) \right) }\text{ if }t>a,
\end{eqnarray*}
which is completely determined because $\overset{.}{x}_{b}\left( t\right)
=x\left( 0\right) \lambda _{d}\left( t\right) e^{\Lambda _{b}\left( t\right)
-\Lambda _{d}\left( t\right) }$ is.

\subsection{Case: $\lambda _{b}\left( a,t\right) $ and $\lambda _{d}\left(
a,t\right) $ independent of time $t$\protect\footnote{%
which amounts to assuming that no modification of the demographical dynamics
ever occurs.}}

Suppose $\lambda _{b}\left( a,t\right) =\lambda _{b}\left( a\right) $ and $%
\lambda _{d}\left( a,t\right) =\lambda _{d}\left( a\right) $ for some given
age-dependent rate functions $\lambda _{b}\left( a\right) ,\lambda
_{d}\left( a\right) .$ Then 
\begin{eqnarray*}
\overline{\Lambda }_{d}\left( a,t\right) &=&\int_{0}^{a\wedge t}\lambda
_{d}\left( a-s\right) ds \\
&=&\Lambda _{d}\left( a\right) -\Lambda _{d}\left( a-t\right) \text{ if }t<a
\\
&=&\Lambda _{d}\left( a\right) \text{ if }t>a.
\end{eqnarray*}
We thus obtain 
\begin{eqnarray*}
\rho \left( a,t\right) &=&\rho _{0}\left( a-t\right) e^{-\left( \Lambda
_{d}\left( a\right) -\Lambda _{d}\left( a-t\right) \right) }\text{ if }t<a \\
&=&\overset{.}{x}_{b}\left( t-a\right) e^{-\Lambda _{d}\left( a\right) }%
\text{ if }t>a
\end{eqnarray*}
which is completely determined as well once $\overset{.}{x}_{b}\left(
t\right) =\rho \left( 0,t\right) $ is. We have 
\begin{eqnarray*}
\overset{.}{x}_{b}\left( t\right) &=&\int \lambda _{b}\left( a\right) \rho
\left( a,t\right) da \\
&=&\int_{0}^{t}\lambda _{b}\left( a\right) \overset{.}{x}_{b}\left(
t-a\right) e^{-\Lambda _{d}\left( a\right) }da+\int_{t}^{\infty }\lambda
_{b}\left( a\right) \rho _{0}\left( a-t\right) e^{-\left( \Lambda _{d}\left(
a\right) -\Lambda _{d}\left( a-t\right) \right) }da,
\end{eqnarray*}
the sum of two convolution terms. Taking the Laplace transforms, with 
\begin{equation*}
\widehat{\alpha }\left( z\right) :=\int_{0}^{\infty }e^{-za}\lambda
_{b}\left( a\right) e^{-\Lambda _{d}\left( a\right) }da\text{ and }\widehat{%
\beta }\left( z\right) :=\int_{0}^{\infty }e^{-za}\rho _{0}\left( a\right)
e^{\Lambda _{d}\left( a\right) }da,
\end{equation*}
we get the Laplace-transform of $\overset{.}{x}_{b}\left( t\right) $ as 
\begin{equation*}
\widehat{x_{b}^{.}}\left( z\right) =\widehat{x_{b}^{.}}\left( z\right) 
\widehat{\alpha }\left( z\right) +\widehat{\beta }\left( z\right) \widehat{%
\alpha }\left( z\right) ,
\end{equation*}
leading to 
\begin{equation*}
\widehat{x_{b}^{.}}\left( z\right) =\frac{\widehat{\alpha }\left( z\right) 
\widehat{\beta }\left( z\right) }{1-\widehat{\alpha }\left( z\right) }.
\end{equation*}
Inverting the Laplace transform yields $\overset{.}{x}_{b}\left( t\right) $.
The solution of (\ref{eq1}) is completely determined from (\ref{eq2}).

We also have

\begin{equation*}
x\left( t\right) =\int \rho \left( a,t\right) da=\int_{0}^{t}\overset{.}{x}%
_{b}\left( t-a\right) e^{-\Lambda _{d}\left( a\right) }da+\int_{t}^{\infty
}\rho _{0}\left( a-t\right) e^{-\left( \Lambda _{d}\left( a\right) -\Lambda
_{d}\left( a-t\right) \right) }da
\end{equation*}
the sum of two convolution terms. Taking the Laplace transforms, with 
\begin{equation*}
\widehat{\alpha }_{0}\left( z\right) :=\int_{0}^{\infty }e^{-za}e^{-\Lambda
_{d}\left( a\right) }da\text{,}
\end{equation*}
we get 
\begin{equation*}
\widehat{x}\left( z\right) =\widehat{x_{b}^{.}}\left( z\right) \widehat{%
\alpha }_{0}\left( z\right) +\widehat{\beta }\left( z\right) \widehat{\alpha 
}_{0}\left( z\right) .
\end{equation*}
So 
\begin{equation}
\widehat{x}\left( z\right) =\frac{\widehat{\alpha }_{0}\left( z\right) 
\widehat{\beta }\left( z\right) }{1-\widehat{\alpha }\left( z\right) }=\frac{%
\widehat{\alpha }_{0}\left( z\right) }{\widehat{\alpha }\left( z\right) }%
\widehat{x_{b}^{.}}\left( z\right) .  \label{sollst}
\end{equation}
Inverting this Laplace transform yields $x\left( t\right) $.

\subsection{The complete case}

This is the realistic case when both $\lambda _{b}\left( a,t\right) $ and $%
\lambda _{d}\left( a,t\right) $ fully depend on $\left( a,t\right) .$
Looking at

\begin{equation*}
\overset{.}{x}\left( t\right) =\int \lambda _{b}\left( a,t\right) \rho
\left( a,t\right) da-\int \lambda _{d}\left( a,t\right) \rho \left(
a,t\right) da=:\text{ }\overset{.}{x}_{b}\left( t\right) -\overset{.}{x}%
_{d}\left( t\right)
\end{equation*}
suggests to introduce the mean birth and death rates (averaging over age) 
\begin{equation*}
\lambda _{b}^{*}\left( t\right) =\frac{\int \lambda _{b}\left( a,t\right)
\rho \left( a,t\right) da}{\int \rho \left( a,t\right) da}\text{ and }%
\lambda _{d}^{*}\left( t\right) =\frac{\int \lambda _{d}\left( a,t\right)
\rho \left( a,t\right) da}{\int \rho \left( a,t\right) da}.
\end{equation*}
Then, with $x\left( t\right) =\int \rho \left( a,t\right) da$%
\begin{equation*}
\overset{.}{x}\left( t\right) =\left( \lambda _{b}^{*}\left( t\right)
-\lambda _{d}^{*}\left( t\right) \right) x\left( t\right) =:\text{ }\overset{%
.}{x}_{b}\left( t\right) -\overset{.}{x}_{d}\left( t\right) .
\end{equation*}
So 
\begin{equation*}
x\left( t\right) =x\left( 0\right) e^{\Lambda _{b}^{*}\left( t\right)
-\Lambda _{d}^{*}\left( t\right) },
\end{equation*}
where $\Lambda _{b}^{*},\Lambda _{d}^{*}$ are the primitives of $\lambda
_{b}^{*},\lambda _{d}^{*}$ and \newline
\begin{equation*}
\overset{.}{x}_{b}\left( t\right) =x\left( 0\right) \lambda _{b}^{*}\left(
t\right) e^{\Lambda _{b}^{*}\left( t\right) -\Lambda _{d}^{*}\left( t\right)
}\text{, }\overset{.}{x}_{d}\left( t\right) =x\left( 0\right) \lambda
_{d}^{*}\left( t\right) e^{\Lambda _{b}^{*}\left( t\right) -\Lambda
_{d}^{*}\left( t\right) }.
\end{equation*}
\newline

We are back to a case where $\lambda _{b}^{*},\lambda _{d}^{*}$ are
independent of age although we know that there is an age dependency for both 
$\lambda _{b}\left( a,t\right) $ and $\lambda _{d}\left( a,t\right) $.

Suppose we are given $\left( \lambda _{b}^{*}\left( t\right) ,\lambda
_{d}^{*}\left( t\right) \right) $ in the first place$.$ Then $\overset{.}{x}%
_{b}\left( t\right) $, $\overset{.}{x}_{d}\left( t\right) $ and $x\left(
t\right) $ are known from Subsection $4.1$. Furthermore $\rho \left(
a,t\right) =\rho ^{*}\left( a,t\right) $ with 
\begin{eqnarray*}
\rho ^{*}\left( a,t\right) &=&\rho _{0}\left( a-t\right) e^{-\Lambda
_{d}^{*}\left( t\right) }\text{ if }t<a \\
&=&\overset{.}{x}_{b}\left( t-a\right) e^{-\left( \Lambda _{d}^{*}\left(
t\right) -\Lambda _{b}^{*}\left( t-a\right) \right) }\text{ if }t>a
\end{eqnarray*}
is the known age-dependent density.

So we need to solve the following inverse problem: find the unknown
functions $\lambda _{b}\left( a,t\right) $ and $\lambda _{d}\left(
a,t\right) $ such that 
\begin{equation*}
\int \lambda _{b}\left( a,t\right) \rho ^{*}\left( a,t\right) da=\text{ }%
\overset{.}{x}_{b}\left( t\right) \text{ and }\int \lambda _{d}\left(
a,t\right) \rho ^{*}\left( a,t\right) da=\text{ }\overset{.}{x}_{d}\left(
t\right)
\end{equation*}
which is a linear problem (not in the convolution class). Then $\left(
x_{b}\left( t\right) ,x_{d}\left( t\right) ,x\left( t\right) \right) $ will
constitute a solution consistent with (\ref{eq1}) and the rates $\lambda
_{b}\left( a,t\right) ,$ $\lambda _{d}\left( a,t\right) .$\newline

We briefly sketch another point of view:

Suppose we were given $\lambda _{d}\left( a,t\right) $ and $\overset{.}{x}%
_{b}\left( t\right) =\rho \left( 0,t\right) $ in the first place$.$

Then $\rho \left( a,t\right) $ is completely known from Equations (\ref
{eq2.0}) and (\ref{eq2}), together with $x\left( t\right) =\int \rho \left(
a,t\right) da$ and $\overset{.}{x}_{d}\left( t\right) =\int \lambda
_{d}\left( a,t\right) \rho \left( a,t\right) da$. Also $x_{b}\left( t\right) 
$ is known from $x_{b}\left( 0\right) =x\left( 0\right) =\int \rho \left(
a,0\right) da$.

Now $\overset{.}{x}_{b}\left( t\right) =\int \lambda _{b}\left( a,t\right)
\rho \left( a,t\right) da$ where now $\lambda _{b}\left( a,t\right) $ has to
be determined. We thus need to solve the inverse problem 
\begin{equation*}
\overset{.}{x}_{b}\left( t\right) =\int_{0}^{t}\lambda _{b}\left( a,t\right) 
\overset{.}{x}_{b}\left( t-a\right) e^{-\overline{\Lambda }^{d}\left(
a,t\right) }da+\int_{t}^{\infty }\lambda _{b}\left( a,t\right) \rho
_{0}\left( a-t\right) e^{-\overline{\Lambda }^{d}\left( a,t\right) }da,
\end{equation*}
where the unknown function is $\lambda _{b}\left( a,t\right) .$ Suppose this
inverse problem was solved. Then $\left( x_{b}\left( t\right) ,x_{d}\left(
t\right) ,x\left( t\right) \right) $ constitute a solution consistent with (%
\ref{eq1}) and the rates $\lambda _{b}\left( a,t\right) ,$ $\lambda
_{d}\left( a,t\right) .$\newline

\textbf{Acknowledgments:} T.H. acknowledges partial support from the labex
MME-DII (Mod\`{e}les Math\'{e}matiques et \'{E}conomiques de la Dynamique,
de l' Incertitude et des Int\'{e}ractions), ANR11-LBX-0023-01.

\end{document}